\documentclass{IEEEtran}
\let\theoremstyle\relax

\usepackage{amsthm}
\theoremstyle{plain}
\newtheorem{theorem}{Theorem}[section]
\newtheorem{lemma}[theorem]{Lemma}
\newtheorem{corollary}[theorem]{Corollary}
\newtheorem{proposition}[theorem]{Proposition}
\theoremstyle{definition}
\newtheorem{definition}[theorem]{Definition}
\theoremstyle{remark}
\newtheorem*{remark}{Remark}
\usepackage{cite}
\usepackage{amsmath,amssymb,amsfonts}
\usepackage{graphicx}
\usepackage{subfigure}
\usepackage{color}

\usepackage{multirow}
\usepackage[ruled,vlined]{algorithm2e}
\usepackage{cuted}

\makeatletter

\newcommand{\Rmnum}[1]{\expandafter\@slowromancap\romannumeral #1@}

\begin{document}
\title{Truncated Multidimensional Trigonometric Moment Problem: A Choice of Bases and the Unique Solution}
\author{Guangyu Wu, \IEEEmembership{Member, IEEE}, Anders Lindquist, \IEEEmembership{Life Fellow, IEEE}
% \thanks{This paragraph of the first footnote will contain the date on 
% which you submitted your paper for review. It will also contain support 
% information, including sponsor and financial support acknowledgment. For 
% example, ``This work was supported in part by the U.S. Department of 
% Commerce under Grant BS123456.'' }
\thanks{This work was completed while Guangyu Wu was with Shanghai Jiao Tong University, Shanghai, China (e-mail: gwu1@alumni.nd.edu).}
\thanks{Anders Lindquist is with School of Artificial Intelligence, Anhui University, Hefei, China (e-mail: alq@kth.se).}}

\maketitle

\begin{abstract}
In this paper, we resolve the Truncated Multidimensional Trigonometric Moment Problem (TMTMP) from a system and signal processing perspective, which serves as the foundation for the Multidimensional Rational Covariance Extension Problem (RCEP). While standard mathematical TMTMPs focus on the existence of atomic measure solutions, system identification requires analytic rational solutions with positive polynomial coefficients. To overcome the long-standing challenge of characterizing the positive feasible domain under general bases, we propose a novel choice of basis functions and a corresponding estimation scheme via convex optimization. We establish an explicit condition to guarantee the positiveness of the spectral estimate. Crucially, the map from the estimate parameters to the trigonometric moments is proved to be a diffeomorphism, ensuring the existence and uniqueness of the solution. Furthermore, we comprehensively prove the statistical properties of the estimator, including its consistency, (asymptotic) unbiasedness, convergence rate, and efficiency. The proposed framework is applied to a system identification task, where simulations validate its effectiveness.
\end{abstract}


\section{Introduction}
\label{sec1}
In this paper we consider the TMTMP in the scope of system and signal processing. The mathematical TMTMP involves determining whether there exists a non-negative measure on the multidimensional torus that matches a given finite set of trigonometric moments. Specifically, given a finite set of complex numbers corresponding to these moments, the challenge is to find a measure such that the integrals of trigonometric polynomial basis functions over this measure equal the given moments. The core goal of the TMTMP is to reconstruct or approximate functions based on partial frequency information. The TMTMP is a cornerstone in mathematical analysis and applied mathematics, highlighting the intricate relationship between algebraic properties and geometric interpretations of measures and moments.

In mathematical TMTMP, the primary concern is often the existence of solutions to the problem. Consequently, solutions in the form of atomic measures are commonly adopted \cite{schmudgen2017moment}. However, in the domain of systems and signal processing, a particular interest lies in a type of TMTMP where the spectral density is parameterized as a rational function with both the numerator and the denominator being positive trigonometric polynomials. From a systems perspective, this form of spectral density enables the modeling of system dynamics as a multidimensional ARMA. Real engineering applications in diverse fields, such as image processing and spectral estimation in radar, sonar, and medical imaging, rely heavily on the development of the theory of this problem. It is also a fundamental theoretical issue in digital twins, which is one of the hottest areas in engineering applications. Spatiotemporal processes are prevalent in digital twins. TMTMP can play a critical role in collecting, integrating, and analyzing data over both space and time to create accurate, dynamic models of real-world systems. This facilitates more efficient management, optimization, and decision-making across various domains. Additionally, data can be insufficient in numerous engineering scenarios. By identifying and modeling the stochastic processes in these scenarios with limited data through TMTMP, it is possible to generate synthetic data samples for training algorithms in digital twin systems, such as reinforcement learning. In conclusion, TMTMP is of great significance in both theory and applications. In the literature, a well-known treatment of the single-dimensional truncated trigonometric moment problem is proposed by the renowned Byrnes-Georgiou-Lindquist school and their academic descendants and collaborators, in which a convex optimization scheme based on minimizing the Kullback-Leibler distance between a reference spectral density and the estimate is utilized \cite{byrnes1998convex, byrnes2001finite, byrnes2002identifiability}. In this paper, we also adopt this scheme of minimizing the functional for the multivariate case. 

We first give a formal definition of the TMTMP considered in this paper, partly following some problem settings in \cite{ringh2016multidimensional}. Given a set of complex numbers $c_{\boldsymbol{k}}, \boldsymbol{k} \in \Omega$, $\boldsymbol{k}:=\left(k_1, \ldots, k_d\right)$ being a vector-valued index belonging to a specified index set $\Omega \subset \mathbb{Z}^d$, which is a finite index set satisfying $\boldsymbol{0} \in \Omega$ and $-\Omega=\Omega$, find a nonnegative bounded measure $d \mu$ such that
\begin{equation}
c_{\boldsymbol{k}}=\int_{\mathbb{T}^d} e^{i(\boldsymbol{k}, \boldsymbol{\theta})} d \mu(\boldsymbol{\theta}) \quad \text { for all } \boldsymbol{k} \in \Omega,
\label{MomentConditions}
\end{equation}
where $\mathbb{T}:=(-\pi, \pi], \boldsymbol{\theta}:=\left(\theta_1, \ldots, \theta_d\right) \in \mathbb{T}^d$, and $(\boldsymbol{k}, \boldsymbol{\theta}):=\sum_{j=1}^d k_j \theta_j$ is the scalar product defined on $\mathbb{R}^d \times \mathbb{R}^d$. Let $e^{i \boldsymbol{\theta}}:=\left(e^{i \theta_1}, \ldots, e^{i \theta_d}\right)$. Moreover, the trigonometric moments satisfy the symmetry condition $c_{-\boldsymbol{k}}=\bar{c}_{\boldsymbol{k}}$. By the Lebesgue decomposition \cite{Rudin_1988}, the measure $d \mu$ can be decomposed in a unique fashion as
$$
d \mu(\boldsymbol{\theta})=\Phi\left(e^{i \boldsymbol{\theta}}\right) d m(\boldsymbol{\theta})+d \breve{\mu}(\boldsymbol{\theta})
$$
with an absolutely continuous part $\Phi d m$ with spectral density $\Phi$ and Lebesgue measure
$$
d m(\boldsymbol{\theta}):=(1 / 2 \pi)^d \prod_{j=1}^d d \theta_j
$$
and a singular part $d \breve{\mu}$ containing, e.g., spectral lines. The measure $d \breve{\mu}$ is an atomic one which is supported on finitely many discrete points. Denote the space of strictly positive continuous functions defined on $\mathbb{T}^{d}$ as $\mathrm{C}^{0}_{+}(\mathbb{T}^{d})$. For the TMTMP, the goal is to parameterize and estimate the finite $\Phi \in \mathrm{C}^{0}_{+}(\mathbb{T}^{d})$ subject to the moment constraints \eqref{MomentConditions}. 

Here we emphasize that the problem we treat in this paper is significantly distinguished from another important research direction of ``multi-dimensional trigonometric moment problem'', where the domain of the estimated spectral density function is $\mathbb{T}$ and the range is a matrix of which each element has range $\mathrm{C}^{0}_{+}(\mathbb{T})$. The basis functions in these problems are chosen as the matricial transfer functions. Representative works on this topic include \cite{georgiou2005solution, georgiou2006relative, ferrante2008hellinger, ferrante2012time, zorzi2013new, zorzi2014multivariate, zorzi2014rational, georgiou2017likelihood, zhu2024statistical}. However, in the problem we treat in this paper, we consider the TMTMP in its original mathematical setting \cite{bakonyi2011matrix}.

We follow the treatment in \cite{ringh2016multidimensional, karlsson2016multidimensional}, where $\Phi$ is obtained by minimizing a generalized entropy functional
\begin{equation}
\mathbb{I}_{P}\left( d\mu \right) = \int_{\mathbb{T}^{d}} P\left(e^{i \boldsymbol{\theta}}\right) \log \frac{P\left(e^{i \boldsymbol{\theta}}\right)}{\Phi\left(e^{i \boldsymbol{\theta}}\right)}d m\left(  \boldsymbol{\theta} \right)
\label{Entropy}
\end{equation}
subject to the moment conditions \eqref{MomentConditions}. It is indeed the Kullback-Leibler distance between the reference density $P\left(e^{i \boldsymbol{\theta}}\right) \in \mathrm{C}^{0}_{+}(\mathbb{T}^{d})$ and the density estimate $\Phi \left(e^{i \boldsymbol{\theta}}\right)$, which we write as $\mathbb{KL}(P\|\Phi)$. In the following sections of this paper, the TMTMP is referred to this minimization problem unless specified otherwise. We note that this optimization problem is infinite-dimensional with a finite number of constraints.

There have been numerous works on this problem \cite{lasserre2009moments, karlsson2016multidimensional, ringh2016multidimensional, ringh2018multidimensional, hao2018equation, zhu2021m2, zhu2021m, zhu2023well, zhu2023weaker}. Quite a great portion of the previous results assume the solution to have a rational form, e.g.
\begin{equation}
\Phi\left(e^{i \boldsymbol{\theta}}\right)=\frac{P\left(e^{i \boldsymbol{\theta}}\right)}{Q\left(e^{i \boldsymbol{\theta}}\right)},
\label{RationalForm}
\end{equation}
where the positive numerator polynomial $P\left(e^{i \boldsymbol{\theta}}\right)$ is given prior. By this assumption, the TMTMP is called a multidimensional rational covariance extension problem. However, this treatment has several weaknesses. First, under the setting of the rational covariance extension problem, $\Phi\left(e^{i \boldsymbol{\theta}}\right)$ in \eqref{RationalForm} is not necessarily optimal. It is possibly suboptimal since the form of solution in \eqref{RationalForm} is assumed \textit{a priori}. Second, in the literature, as to ensure that $\Phi\left(e^{i \boldsymbol{\theta}}\right)$ in the form of \eqref{RationalForm} is nonnegative, it is always assumed that $P\left(e^{i \boldsymbol{\theta}}\right), Q\left(e^{i \boldsymbol{\theta}}\right) \in \overline{\mathfrak{P}}_{+} \backslash \{ 0 \}$, where $\mathfrak{P}_{+}$ is the convex cone of positive trigonometric polynomials
\begin{equation}
Q\left(e^{i \boldsymbol{\theta}}\right)=\sum_{\boldsymbol{k} \in \Omega} \lambda_{\boldsymbol{k}} e^{-i(\boldsymbol{k}, \boldsymbol{\theta})}
\label{Qeitheta}
\end{equation}
that are positive for all $\boldsymbol{\theta} \in \mathbb{T}^d$, and $\overline{\mathfrak{P}}_{+}$ is its closure. However, the choice of $\Omega$ is a very challenging task, since useful descriptions of strictly positive multidimensional trigonometric polynomials are missing \cite{schmudgen2017moment}. Third, given a specific choice of $\Omega$, namely a set of basis functions, the feasible domains of the parameters of $P\left(e^{i \boldsymbol{\theta}}\right), Q\left(e^{i \boldsymbol{\theta}}\right)$, for example $\left\{ \lambda_{\boldsymbol{k}}: \boldsymbol{k} \in \Omega \right\}$ for $Q\left(e^{i \boldsymbol{\theta}}\right)$, are implicit and difficult. Even though the feasible domains are in general convex, however, it might be a null set given some choices of $\Omega \subset \mathbb{Z}^d$ \cite{karlsson2016multidimensional}.

In this paper, we present a novel approach to address the TMTMP, effectively resolving the three issues in the preceding results. The contribution of our proposed scheme falls within the following three aspects.

(1) A choice of bases. We propose a choice of $\Omega$, and prove that given this choice of basis functions, $\Phi\left(e^{i \boldsymbol{\theta}}\right)$ of the form \eqref{RationalForm} is the optimal solution to the TMTMP. It is to say that given our proposed $\Omega$, the TMTMP is essentially the multidimensional rational covariance extension problem, without assuming \textit{a priori} that $\Phi\left( e^{i \boldsymbol{\theta}} \right)$ has the form in \eqref{RationalForm}.

(2) Positiveness. We propose the condition for our parameterization of $Q\left( e^{i\boldsymbol{\theta}} \right)$ to be positive on $\mathbb{T}^{d}$. The condition, namely the feasible domain of the parameters $\left\{ \lambda_{\boldsymbol{k}}: \boldsymbol{k} \in \Omega \right\}$ in $Q\left(e^{i \boldsymbol{\theta}}\right)$, is explicit.

(3) Convex optimization. We prove that the feasible domain of $\left\{ \lambda_{\boldsymbol{k}}: \boldsymbol{k} \in \Omega \right\}$ is a convex set. Moreover, we prove that the cost function is a convex one. Hence maximizing \eqref{Entropy} subject to \eqref{MomentConditions} is a convex optimization problem. We also prove the existence and uniqueness of solution to this optimization problem.

The paper is organized as follows. In Section 2, we propose a choice of basis functions and a corresponding convex optimization scheme. However, unlike the one-dimensional trigonometric problem, the positiveness of the solution of the problem is not ensured. We then put forward a positive parametrization of the multidimensional spectral density subsequently in Section 3. In Section 4, we prove the existence and uniqueness of the solution to the proposed convex optimization problem. We consider using the proposed parametrization for estimation in Section 5, where two types of statistics are given. A necessary and sufficient condition for the existence of solution of the convex optimization using the statistics is proposed. Furthermore, consistency, (asymptotic) unbiasedness, convergence rate and efficiency under a mild assumption are proved. An algorithm for spectral density estimation is put forward. Simulation results are given in Section 6, including generating samples from the true spectral density and estimating the spectral density. The results of the proposed algorithm are compared to minimizing the L2 norm of error, which reveals the advantage of the proposed algorithm. 

\section{A choice of basis functions and the convex optimization scheme}
\label{sec2}

In most previous results on the multidimensional trigonometric moment problems, for example \cite{karlsson2016multidimensional, ringh2016multidimensional}, the basis functions of the parametrization for the TMTMP are not determined explicitly, i.e., a choice of the elements in the set $\Omega$ is not given explicitly. The TMTMP is defined as a very general problem in these works, where the number and forms of the basis functions are quite flexible. However, we note that the general formulation of the TMTMP, e.g. minimizing $\mathbb{I}_{P}\left( d\mu \right)$ subject to the moment constraints in Section \ref{sec1}, is not a well-defined one. A specific choice of the basis functions doesn't necessarily yield a positive parametrization of the spectral density function, which means that the choice of the basis functions is not arbitrary. However, to the best of our knowledge, there has not been a result on the choice of a specific set of basis functions for the TMTMP and the corresponding condition for a parametrization to be positive. In this section and the next one, we will propose a parametrization of the TMTMP with a selected set of basis functions, together with the positiveness proof.

A common choice of $\Omega$ is
\begin{equation}
\Omega := \left\{ \left(k_{1}, \cdots, k_{d} \right) \mid |k_{i}| \leqslant n, n \in \mathbb{N}, i = 1, \cdots, d \right\},
\label{OmegaSet}
\end{equation}
where $\mathbb{N}$ denotes the set of natural numbers, $n$ is the prescribed highest order of the trigonometric moments \cite{ringh2018multidimensional}. In the TMTMP we treat, we also consider the $\Omega$ in \eqref{OmegaSet}. Now the problem comes to proving that the parametrization using the selected basis functions is positive under some explicit conditions. 

Given the choice of basis functions, we are then provided with the trigonometric moments $c_{\boldsymbol{k}}$ corresponding to the basis functions, which have the form in \eqref{MomentConditions}. However, since we don't assume the analytic form of the parametrization of the spectral density estimate, it is not feasible for us to obtain the estimate directly by the trigonometric moments. We then follow the treatment in \cite{karlsson2016multidimensional, ringh2016multidimensional} by minimizing the generalized entropy functional \eqref{Entropy}.

We note that the functional \eqref{Entropy}, where $\Phi$ is the absolutely continuous part of $d \mu$, is convex, but not strictly convex, because the singular part of the measure does not alter the value. Since the trigonometric constraints \eqref{MomentConditions} are linear, this optimization problem is a convex one. However, since the optimization problem is infinite-dimensional, it is more convenient to consider the dual problem. The dual functional of minimizing a generalized entropy functional \eqref{Entropy} given the trigonometric moment conditions \eqref{MomentConditions} can be written as
\begin{equation}
\begin{aligned}
& \mathcal{L}\left(\Phi, \left\{ \lambda_{\boldsymbol{k}} : \boldsymbol{k} \in \Omega \right\} \right)\\
= & \int_{\mathbb{T}^d} P\left(e^{i \boldsymbol{\theta}}\right) \log \frac{P\left(e^{i \boldsymbol{\theta}}\right)}{\Phi\left(e^{i \boldsymbol{\theta}}\right)} d m(\boldsymbol{\theta})\\
+ & \sum_{\boldsymbol{k} \in \Omega} \lambda_{\boldsymbol{k}} \left( \int_{\mathbb{T}^d} e^{i(\boldsymbol{k}, \boldsymbol{\theta})} \Phi\left(e^{i \boldsymbol{\theta}}\right) d m(\boldsymbol{\theta}) - c_{\boldsymbol{k}}  \right),
\end{aligned}
\label{DualFunctional}
\end{equation}
where $\left\{ \lambda_{\boldsymbol{k}} : \boldsymbol{k} \in \Omega \right\}$ are the Lagrangian multipliers. Now the problem is to maximize the dual functional
\begin{equation}
\left\{ \lambda_{\boldsymbol{k}} : \boldsymbol{k} \in \Omega \right\} \mapsto \inf_{\Phi \in \mathrm{C}^{0}_{+}(\mathbb{T}^{d})}  \mathcal{L}\left(\Phi, \left\{ \lambda_{\boldsymbol{k}} : \boldsymbol{k} \in \Omega \right\} \right).
\label{DualFunc}
\end{equation}
Clearly $\Phi \mapsto \mathcal{L}\left(\Phi, \left\{ \lambda_{\boldsymbol{k}} : \boldsymbol{k} \in \Omega \right\} \right)$ is strictly convex. The primal optimization problem \eqref{Entropy} under the trigonometric moment constraints is a problem in infinite dimensions but with a finite number of constraints. The dual to this problem will then have a finite number of variables but an infinite number of constraints. To be able to determine the RHS of \eqref{DualFunc}, we need to find a feasible $\Phi$, for which the directional derivative $\delta \mathcal{L}\left(\Phi, \left\{ \lambda_{\boldsymbol{k}} : \boldsymbol{k} \in \Omega \right\}; \delta \Phi \right)=0$ for all relevant $\delta \Phi$. This further restricts the choice of $\left\{ \lambda_{\boldsymbol{k}} : \boldsymbol{k} \in \Omega \right\}$. The directional derivative reads
$$
\begin{aligned}
& \delta \mathcal{L}\left(\Phi, \left\{ \lambda_{\boldsymbol{k}} : \boldsymbol{k} \in \Omega \right\}; \delta \Phi \right)\\
= & \int_{\mathbb{T}^d} \left( - \frac{P\left(e^{i \boldsymbol{\theta}}\right)}{ \Phi\left(e^{i \boldsymbol{\theta}}\right)} + \sum_{\boldsymbol{k} \in \Omega} \lambda_{\boldsymbol{k}} e^{i(\boldsymbol{k}, \boldsymbol{\theta})} \right) \delta \Phi\left(e^{i \boldsymbol{\theta}}\right) d m(\boldsymbol{\theta}).
\end{aligned}
$$

It has to be zero at a minimum for all variations $\delta \Phi$. We note that this can be achieved only if
\begin{equation}
\Phi\left(e^{i \boldsymbol{\theta}}\right) = \frac{P\left(e^{i \boldsymbol{\theta}}\right)}{\sum_{\boldsymbol{k} \in \Omega} \lambda_{\boldsymbol{k}} e^{i(\boldsymbol{k}, \boldsymbol{\theta})}}
\label{PhiForm}
\end{equation}
for all $\boldsymbol{\theta} \in \mathbb{T}^{d}$. We may note that the density estimate $\Phi(e^{i \boldsymbol{\theta}})$ has a similar form as that in \eqref{RationalForm}. However, $Q(e^{i \boldsymbol{\theta}}) \in \overline{\mathfrak{P}}_{+} \backslash \{ 0 \}$ in the rational form of spectral density estimate \eqref{RationalForm}. Therefore, there still remain two key problems. The first is to derive an explicit condition by which we have
$$
\sum_{\boldsymbol{k} \in \Omega} \lambda_{\boldsymbol{k}} e^{i(\boldsymbol{k}, \boldsymbol{\theta})} \in \overline{\mathfrak{P}}_{+} \backslash \{ 0 \}.
$$
The second is to prove that given the proposed condition, the optimization problem to estimate the spectral density has a solution, and the solution is unique. In the following parts of the paper, we will first consider the positiveness condition of $\Phi\left(e^{i \boldsymbol{\theta}}\right)$, which guarantees the positiveness of $\Phi\left(e^{i \boldsymbol{\theta}}\right)$ in \eqref{PhiForm}. Then the existence and uniqueness of solution to the optimization problem will be given.

We first review some results from the one-dimensional trigonometric moment problem, which has the same problem setting as that in Section \ref{sec1}, with $d = 1$. In \cite{georgiou2003kullback}, the basis functions are chosen as $\Omega = \left\{ 0, 1, \cdots, n \right\}$. Denote
$$
K(e^{i\theta}) := \begin{bmatrix}
1 & e^{i\theta} & \cdots & e^{i(n-1)\theta} & e^{in\theta}
\end{bmatrix}^{\intercal}.
$$
With the selected set of bases, the trigonometric moment conditions can be written in the form of a matrix equation
\begin{equation}
\int_{\mathbb{T}} K(e^{i\theta}) \Phi(e^{i\theta}) K^{H}(e^{i\theta}) d m(\theta)=\mathcal{T}_{1},
\label{MatrixCondition}
\end{equation}
where $M^{H}$ denotes the conjugate transpose of a matrix $M$. The right member of \eqref{MatrixCondition} is the Toeplitz matrix
$$
\mathcal{T}_{1} = \left[\begin{array}{cccc}
c_0 & c_1 & \cdots & c_n \\
\bar{c}_1 & c_0 & \cdots & c_{n-1} \\
\vdots & \vdots & \ddots & \vdots \\
\bar{c}_n & \bar{c}_{n-1} & \cdots & c_0
\end{array}\right],
$$
where $\bar{c}$ denotes the conjugate of a complex number $c$. Moreover, $\mathcal{T}_{1}$ is positive definite. In fact, $\mathcal{T}_{1}$ is in the range of an integral operator
$$
\Gamma_{1}: \Phi(e^{i\theta}) \mapsto \mathcal{T}_{1} = \int_{\mathbb{T}} K(e^{i\theta}) \Phi(e^{i\theta}) K^{H}(e^{i\theta}) d m(\theta).
$$

Define
\begin{equation}
    \mathfrak{L}^{1}_{+}
    := \left\{\Lambda \in \operatorname{range}(\Gamma_{1}) \mid \Lambda \succ 0 \right\}.
\end{equation}

By Theorem 5 in \cite{georgiou2003kullback}, we have that with $d = 1$, the truncated trigonometric moment problem, i.e., minimizing \eqref{Entropy} under the trigonometric moment conditions \eqref{MomentConditions}, has a unique solution in the form of
$$
\Phi\left(e^{i \boldsymbol{\theta}}\right)=\frac{P\left(e^{i \boldsymbol{\theta}}\right)}{K^{H}\left(e^{i \boldsymbol{\theta}}\right) \Lambda K\left(e^{i \boldsymbol{\theta}}\right)},
$$
where $\Lambda \in \mathfrak{L}^{1}_{+}$ is the Lagrangian matrix. By constraining $\Lambda$ to fall within $\mathfrak{L}^{1}_{+}$, we have $\Phi\left(e^{i \boldsymbol{\theta}}\right) > 0$.

This exquisite treatment by writing the trigonometric moments in the form of a matrix equation utilizes the property of positive definiteness, which provides an explicit condition for the estimate to be strictly positive. However, it is only valid for $d =1$. For the TMTMP, there doesn't exist a positive parametrization of the spectral density estimate, with an explicit condition of positiveness, to our knowledge. In the next section, we will propose such a positive parametrization with a corresponding positiveness condition for the TMTMP.

\section{A positive multivariate spectral parametrization}
\label{sec3}

Denote the Kronecker product as $\otimes$. We write
$$
K(e^{i\theta_{k}}) := \begin{bmatrix}
1 & e^{i\theta_{k}} & \cdots & e^{in\theta_{k}}
\end{bmatrix}^{\intercal}, \quad k = 1, \cdots, d.
$$
Then the bases \eqref{OmegaSet} we chose in Section \ref{sec2} can be represented by a series of Kronecker products
\begin{equation}
\begin{aligned}
    & K(e^{i\boldsymbol{\theta}})\\
= & K\left(e^{i\theta_{1}}\right) \otimes K\left(e^{i\theta_{2}}\right) \otimes \cdots \otimes K\left(e^{i\theta_{d}}\right).
\end{aligned}
\end{equation}
Moreover, we note that the spectral density estimate \eqref{PhiForm} can be written in the form
\begin{equation}
    \Phi\left(e^{i \boldsymbol{\theta}}\right) = \frac{P\left(e^{i \boldsymbol{\theta}}\right)}{K^{H}(e^{i\boldsymbol{\theta}}) \Lambda K(e^{i\boldsymbol{\theta}})}.
\label{MatrixPhi}
\end{equation}
We note that there are $(n+1)^{2d}$ entries (Lagrangian multipliers) in the matrix $\Lambda$. However, there are $(n+1)^{d}$ trigonometric moment constraints in total (we regard $c_{\boldsymbol{k}}$ and its conjugate $\bar{c}_{\boldsymbol{k}}$ as the same constraint). It is obvious that the latter one is always less than the former one with $n, k \in \mathbb{N}_0 \setminus \{0\}$. Therefore, some entries of $\Lambda$ are identical and correspond to the same trigonometric moment condition. Denote the $i_{\text{th}}$ element of $K(e^{i\boldsymbol{\theta}})$ as $K_{i}$, and the entry of $\Lambda$ in row $i$ and column $j$ as $\Lambda_{i, j}$. Let the set
$$
\mathcal{N}_{\boldsymbol{k}}:=\left\{\Lambda_{i, j} \mid K_i \cdot \bar{K}_{j}=e^{ik_{1}\theta_{1}}e^{ik_{2}\theta_{2}} \cdots e^{ik_{d}\theta_{d}}\right\} .
$$
We note that all $\Lambda_{i, j} \in \mathcal{N}_{\boldsymbol{k}}$ correspond to the same moment constraint
\begin{equation}
c_{\boldsymbol{k}}=\int_{\mathbb{T}^d} e^{ik_{1}\theta_{1}} e^{ik_{2}\theta_{2}} \cdots e^{ik_{d}\theta_{d}} \Phi\left(e^{i \boldsymbol{\theta}}\right) d m(\boldsymbol{\theta}),
\label{ck}
\end{equation}
hence are identical. Denote the cardinality of $\mathcal{N}_{\boldsymbol{k}}$ as $\operatorname{card}(\mathcal{N}_{\boldsymbol{k}})$. For each $\Lambda_{i, j} \in \mathcal{N}_{\boldsymbol{k}}$, we have \footnote{It is not necessary to calculate $\operatorname{card}(\mathcal{N}_{\boldsymbol{k}})$ for the subsequent optimization scheme. Equation \eqref{LambdaBreve} is provided solely to demonstrate the relationship between \eqref{PhiForm} and \eqref{MatrixPhi}. It suffices to ensure that all $\Lambda_{i, j}$ corresponding to the same $\mathcal{N}_{\boldsymbol{k}}$ are identical.}
\begin{equation}
\Lambda_{i, j}=\frac{\lambda_{\boldsymbol{k}}}{\operatorname{card}(\mathcal{N}_{\boldsymbol{k}})}.
\label{LambdaBreve}
\end{equation}
Moreover, we have that $\mathcal{T}_{d}$ is in the range of the linear integral operator
\begin{equation}
\Gamma_{d}: \Phi\left(e^{i \boldsymbol{\theta}}\right) \mapsto \mathcal{T}_{d}=\int_{\mathbb{T}^d} K(e^{i \boldsymbol{\theta}}) \Phi\left(e^{i \boldsymbol{\theta}}\right) K^{H}(e^{i \boldsymbol{\theta}}) d m(\boldsymbol{\theta}).
\label{MomentConstraintsEq}
\end{equation}
We have given a form of parametrization for the spectral density in \eqref{MatrixPhi}. It remains to propose a condition for the parametrization to be strictly positive. To our knowledge, there has not been such a result of positiveness for the TMTMP in the literature. 

Define
\begin{equation}
    \mathfrak{L}^{d}_{+}
    := \left\{\Lambda \in \operatorname{range}(\Gamma_{d}) \mid \Lambda \succ 0 \right\},
\end{equation}
and we have the following lemma.

\begin{lemma}
$\Phi\left(e^{i \boldsymbol{\theta}}\right) > 0$ for all $\boldsymbol{\theta} \in \mathbb{T}^{d}$ only if $\Lambda \in \mathfrak{L}^{d}_{+}$.
\label{lemma31}
\end{lemma}
\begin{proof}
    Since $P\left(e^{i \boldsymbol{\theta}}\right)$ in \eqref{MatrixPhi} is strictly positive, it is equivalent to proving that $ K^{H}(e^{i\boldsymbol{\theta}}) \Lambda K(e^{i\boldsymbol{\theta}}) > 0$ for all $\boldsymbol{\theta} \in \mathbb{T}^{d}$ only if $\Lambda \in \operatorname{range}(\Gamma_{d})$. 

    We have that $\Phi\left(e^{i \boldsymbol{\theta}}\right)$ is strictly positive on $\mathbb{T}^{d}$, and 
    $$
    K(e^{i \boldsymbol{\theta}}) K^{H}(e^{i \boldsymbol{\theta}}) \succ 0, \quad \forall \boldsymbol{\theta} \in \mathbb{T}^{d}.
    $$ 
    Since the sum of positive definite matrices is still positive definite, we have that
    $$
        \int_{\mathbb{T}^d} K(e^{i \boldsymbol{\theta}}) \Phi\left(e^{i \boldsymbol{\theta}}\right) K^{H}(e^{i \boldsymbol{\theta}}) d m(\boldsymbol{\theta}) \succ 0,
    $$
    i.e., all elements within the set $\operatorname{range}(\Gamma_{d})$ are positive definite. Therefore, we have $\Lambda \succ 0$.
    
    By the property of positive definite matrix, see for example \cite{franklin2012matrix}, we have 
    $K^{H}(e^{i\boldsymbol{\theta}}) \Lambda K(e^{i\boldsymbol{\theta}}) > 0$ with $K(e^{i\boldsymbol{\theta}}) \neq \boldsymbol{0}$,
    given $\Lambda \succ 0$. We note that $K(e^{i\boldsymbol{\theta}}) \neq \boldsymbol{0}, \forall  \boldsymbol{\theta} \in \mathbb{T}^{d}$. The proof is then complete.
\end{proof}

Lemma \ref{lemma31} provides a necessary condition for the spectral density estimate being positive, which means that the set $\operatorname{range}(\Gamma_{d})$ is a subset of the set of all feasible $\Lambda$ satisfying $\Phi\left(e^{i \boldsymbol{\theta}}\right) > 0$. In the following section, we will prove that the range of optimal $\Lambda$, which maximizes the dual functional \eqref{DualFunctional}, is essentially the set $\mathfrak{L}^{d}_{+}$. Prior to that, we would first prove that the set $\mathfrak{L}^{d}_{+}$ is convex.

\begin{lemma}
    The set $\mathfrak{L}^{d}_{+}$ is convex.
    \label{Lemma32}
\end{lemma}

\begin{proof}
    Let $\Lambda_{1}, \Lambda_{2}$ be two arbitrary elements of $\mathfrak{L}^{d}_{+}$. It suffices to prove
    $$
        \alpha_{1}\Lambda_{1} + \alpha_{2}\Lambda_{2} \in \mathfrak{L}^{d}_{+}
    $$
    for all pairs of $\alpha_{1}, \alpha_{2} > 0$ satisfying $\alpha_{1}+\alpha_{2}=1$.

    For $\Lambda_{1}, \Lambda_{2} \in \mathfrak{L}^{d}_{+}$, there exist $\Phi_{1}\left(e^{i \boldsymbol{\theta}}\right), \Phi_{2}\left(e^{i \boldsymbol{\theta}}\right) \in \mathrm{C}^{0}_{+}(\mathbb{T}^{d})$ respectively. We can then write
    $$
    \begin{aligned}
        & \alpha_{1}\Lambda_{1} + \alpha_{2}\Lambda_{2}\\
        = & \int_{\mathbb{T}^d} \alpha_{1} K(e^{i \boldsymbol{\theta}}) \Phi_{1}\left(e^{i \boldsymbol{\theta}}\right) K^{H}(e^{i \boldsymbol{\theta}}) d m(\boldsymbol{\theta})\\
        + & \int_{\mathbb{T}^d} \alpha_{2} K(e^{i \boldsymbol{\theta}}) \Phi_{2}\left(e^{i \boldsymbol{\theta}}\right) K^{H}(e^{i \boldsymbol{\theta}}) d m(\boldsymbol{\theta})\\
        = & \int_{\mathbb{T}^d} K(e^{i \boldsymbol{\theta}}) \left(\alpha_{1}\Phi_{1}\left(e^{i \boldsymbol{\theta}}\right)+\alpha_{2}\Phi_{2}\left(e^{i \boldsymbol{\theta}}\right)\right)K^{H}(e^{i \boldsymbol{\theta}}) d m(\boldsymbol{\theta}).
    \end{aligned}
    $$
    We note that $\alpha_{1}\Phi_{1}\left(e^{i \boldsymbol{\theta}}\right)+\alpha_{2}\Phi_{2}\left(e^{i \boldsymbol{\theta}}\right) \in \mathrm{C}^{0}_{+}(\mathbb{T}^{d})$. Therefore, $\alpha_{1}\Lambda_{1} + \alpha_{2}\Lambda_{2} \in \mathfrak{L}^{d}_{+}$, which completes the proof.
\end{proof}

By Lemma \ref{Lemma32}, we have proved that $\mathfrak{L}^{d}_{+}$ is a convex subset of the set of all positive definite matrices. 
\section{Existence and uniqueness of the solution}

In the previous section, we have proposed a multivariate spectral parametrization by power moments. Moreover, the parametrization \eqref{MatrixPhi} is proved to be strictly positive with $\Lambda$ to fall within the set $\mathfrak{L}^{d}_{+}$. However, the optimal solution which maximizes the dual functional \eqref{DualFunctional} has not been proved to fall within $\mathfrak{L}^{d}_{+}$. In this section, we prove the existence and uniqueness of the solution to maximizing the dual functional \eqref{DualFunctional}, with the domain of $\Lambda$ being $\mathfrak{L}^{d}_{+}$. 

By the moment conditions, in the form of a matrix equation \eqref{MomentConstraintsEq}, the Lagrangian \eqref{DualFunctional} can then be adapted to
\begin{equation}
\begin{aligned}
\mathcal{L}(\Phi, \Lambda)
= & \mathbb{I}_{P}\left( d\mu \right)+ \operatorname{tr}\left (\Lambda(\Gamma_{d}(\Phi)-\mathcal{T}_{d}) \right)\\
= & \int_{\mathbb{T}^{d}} P\left(e^{i \boldsymbol{\theta}}\right) \log \frac{P\left(e^{i \boldsymbol{\theta}}\right)}{\Phi\left(e^{i \boldsymbol{\theta}}\right)}d m\left(  \boldsymbol{\theta} \right) - \operatorname{tr}\left( \Lambda \mathcal{T}_{d} \right)\\
+ & \operatorname{tr}\left( \int_{\mathbb{T}^d} \Lambda K(e^{i \boldsymbol{\theta}}) \Phi\left(e^{i \boldsymbol{\theta}}\right) K^{H}(e^{i \boldsymbol{\theta}}) d m(\boldsymbol{\theta})\right).
\end{aligned}
\label{LrhoLambda}
\end{equation}

By substituting \eqref{MatrixPhi} into \eqref{LrhoLambda}, we obtain
$$
\begin{aligned}
& \mathcal{L}(\Phi, \Lambda)\\
= & \int_{\mathbb{T}^{d}} P\left(e^{i \boldsymbol{\theta}}\right) \log \left( K^{H}(e^{i\boldsymbol{\theta}}) \Lambda K(e^{i\boldsymbol{\theta}}) \right)d m\left(  \boldsymbol{\theta} \right) - \operatorname{tr}\left( \Lambda \mathcal{T}_{d} \right)\\
+ & \operatorname{tr} \left( \int_{\mathbb{T}^{d}}\frac{\Lambda K(e^{i \boldsymbol{\theta}}) P\left(e^{i \boldsymbol{\theta}}\right) K^{H}(e^{i \boldsymbol{\theta}})}{K^{H}(e^{i \boldsymbol{\theta}})\Lambda K(e^{i \boldsymbol{\theta}})} d m(\boldsymbol{\theta}) \right)\\
= & \int_{\mathbb{T}^{d}} \left(P(e^{i \boldsymbol{\theta}}) \log \left(K^{H}(e^{i \boldsymbol{\theta}})\Lambda K(e^{i \boldsymbol{\theta}})\right) + P(e^{i \boldsymbol{\theta}}) \right) dm(\boldsymbol{\theta})\\
- & \operatorname{tr}\left( \Lambda \mathcal{T}_{d} \right).
\end{aligned}
$$

Then the problem turns to minimizing the dual functional
\begin{equation}
\begin{aligned}
& \mathbb{J}_{P}(\Lambda)\\
:= & - \int_{\mathbb{T}^{d}} P(e^{i \boldsymbol{\theta}}) \log\left(K^{H}(e^{i \boldsymbol{\theta}})\Lambda K(e^{i \boldsymbol{\theta}})\right) dm(\boldsymbol{\theta})
+ \operatorname{tr}(\Lambda \mathcal{T}_{d}).
\end{aligned}
\label{LossFunc}
\end{equation}

We then have the following lemma.

\begin{lemma}
Define the map $\omega: \mathfrak{L}^{d}_{+} \rightarrow \mathfrak{S}^{d}_{+}$ between $\mathfrak{L}^{d}_{+}$ and $\mathfrak{S}^{d}_{+}:=\{\mathcal{T}_{d} \in \operatorname{range}(\Gamma_{d}) \mid \mathcal{T}_{d} \succ 0\}$. Any stationary point of $\mathbb{J}_{P}(\Lambda)$ must satisfy
\begin{equation}
    \omega(\Lambda)=\mathcal{T}_{d},
    \label{Omega}
\end{equation}
where we have defined the map
$$
\omega: \; \Lambda \mapsto \int_{\mathbb{T}^{d}} K(e^{i \boldsymbol{\theta}}) \frac{P(e^{i \boldsymbol{\theta}})}{Q \left( e^{i \boldsymbol{\theta}}, \Lambda \right)} K^{H}(e^{i \boldsymbol{\theta}})d m(\boldsymbol{\theta})
\label{omega}
$$
with $Q \left( e^{i \boldsymbol{\theta}}, \Lambda \right):=K^{H}(e^{i\boldsymbol{\theta}}) \Lambda K(e^{i\boldsymbol{\theta}})$. 
\label{Lemma42}
\end{lemma}

\begin{proof}
The directional derivative of $Q \left( e^{i \boldsymbol{\theta}}, \Lambda \right)$ with respect to $\Lambda$ reads
\begin{displaymath}
\delta Q(e^{i\boldsymbol{\theta}}, \Lambda ; \delta \Lambda)=K^{H} \delta \Lambda K=\operatorname{tr}\left(\delta \Lambda K K^{H}\right),
\end{displaymath}
by which we obtain the directional derivative
\begin{equation*}
\begin{aligned}
& \delta \mathbb{J}_{P}(\Lambda ; \delta \Lambda)\\
= & - \int_{\mathbb{T}^{d}} \frac{P(e^{i \boldsymbol{\theta}})}{Q \left( e^{i \boldsymbol{\theta}}, \Lambda \right)} \delta Q \left( e^{i \boldsymbol{\theta}}, \Lambda; \delta \Lambda \right) dm(\boldsymbol{\theta}) + \operatorname{tr}\left( \delta \Lambda \mathcal{T} \right)\\
= & \operatorname{tr}\left(\delta \Lambda\left[\mathcal{T}-\int_{\mathbb{T}^{d}} K(e^{i \boldsymbol{\theta}})\frac{P(e^{i \boldsymbol{\theta}})}{Q(e^{i \boldsymbol{\theta}}, \Lambda)} K^{H}(e^{i \boldsymbol{\theta}})dm(\boldsymbol{\theta})\right]\right)
\end{aligned}
\label{FirstOrder}
\end{equation*}
being zero for all $\delta \Lambda \in\operatorname{range}(\Gamma_{d})$ if and only if \eqref{Omega} holds, which completes the proof. 
\end{proof}

To establish the existence and uniqueness of the minimum of $\mathbb{J}_{P}$, it is necessary to demonstrate that the map $\omega: \mathfrak{L}^{d}_{+} \rightarrow \mathfrak{S}^{d}_{+}$ is both injective, thereby establishing uniqueness; and surjective, thereby establishing existence. Consequently, we prove that \eqref{Omega} admits a unique solution, hence implying the existence of a unique minimum for the dual functional $\mathbb{J}_{P}$. We initiate by considering injectivity.

\begin{lemma}
Suppose $\Lambda \in \text{range}(\Gamma_{d})$. Then the map

\begin{equation}
\label{Lambda2G'LambdaG}
\Lambda \mapsto K^{H} \Lambda K
\end{equation}
is injective.
\label{Lemma43}
\end{lemma}

\begin{proof} 
Since $\Lambda \in \text{range}(\Gamma_{d})$,
$$
   \Lambda = \int_{\mathbb{T}^{d}} K(e^{i\boldsymbol{\gamma}}) \psi(e^{i\boldsymbol{\gamma}}) K^{H}(e^{i\boldsymbol{\gamma}})dm(\boldsymbol{\gamma})
$$
for some $\psi \in \mathrm{C}^{0}_{+}(\mathbb{T}^{d})$.
Suppose $K^{H}\Lambda K = 0$. Then we have $\int_{\mathbb{T}^{d}} K^{H}(e^{i\boldsymbol{\theta}})\Lambda K(e^{i\boldsymbol{\theta}})dm(\boldsymbol{\theta}) = 0$, and therefore we have \eqref{longeq1}, displayed at the top of the next page.

\begin{figure*}[t]
\begin{equation}
\begin{aligned}
    & \int_{\mathbb{T}^{d}} K^{H}(e^{i\boldsymbol{\theta}})\Lambda K(e^{i\boldsymbol{\theta}})dm(\boldsymbol{\theta}) \\
    = & \operatorname{tr} \left( \int_{\mathbb{T}^{d}} K^{H}(e^{i\boldsymbol{\theta}}) \int_{\mathbb{T}^{d}} K(e^{i\boldsymbol{\gamma}}) \psi(e^{i\boldsymbol{\gamma}}) K^{H}(e^{i\boldsymbol{\gamma}})dm(\boldsymbol{\gamma})K(e^{i\boldsymbol{\theta}})dm(\boldsymbol{\theta}) \right) \\
    = & \int_{\mathbb{T}^{d}} \int_{\mathbb{T}^{d}} [K^{H}(e^{i\boldsymbol{\theta}})K(e^{i\boldsymbol{\gamma}})]^2\psi(e^{i\boldsymbol{\gamma}}) dm(\boldsymbol{\gamma})dm(\boldsymbol{\theta}) =0.
\end{aligned}
\label{longeq1}
\end{equation}
\hrulefill
\vspace*{4pt}
\end{figure*}

Consequently we have $[K^{H}(e^{i\boldsymbol{\theta}})K(e^{i\boldsymbol{\gamma}})]^2\psi(e^{i\boldsymbol{\gamma}}) = 0$, for all $\boldsymbol{\theta},\boldsymbol{\gamma} \in \mathbb{T}^{(n+1)^{d}}$, which clearly implies that $\psi(e^{i\boldsymbol{\gamma}})=0$, and hence that $\Lambda = 0$.
Therefore, the map \eqref{Lambda2G'LambdaG} is injective, as claimed.
\end{proof}

\begin{theorem}
The map $\omega: \mathfrak{L}^{d}_{+} \mapsto \mathfrak{S}^{d}_{+}$ is a homeomorphism. 
\label{Theorem43}
\end{theorem}
\begin{proof}
We first prove that $\omega: \mathfrak{L}^{d}_{+} \mapsto \mathfrak{S}^{d}_{+}$ is injective. Suppose that $\omega(\Lambda_1)=\omega(\Lambda_2)$ for some $\Lambda_1$ and $\Lambda_2$ in $\mathfrak{L}^{d}_+$. We need to show that $\Lambda_1=\Lambda_2$. We note that
$$
\begin{aligned}
\Delta \omega = & \omega(\Lambda_1)-\omega(\Lambda_2)\\
= & \int_{\mathbb{T}^{d}}K(e^{i\boldsymbol{\theta}})K^{H}(e^{i\boldsymbol{\theta}})\frac{P(e^{i\boldsymbol{\theta}})}{Q_{1}Q_{2}}(Q_{2}-Q_{1})dm(\boldsymbol{\theta})\\
= & 0,
\end{aligned}
$$
where $Q_{1}=K^{H}\Lambda_{1}K^{H}$ and $Q_{2}=K^{H}\Lambda_{2}K^{H}$. The element of $\Delta \omega$ in the first row and the first column reads
\begin{equation}
\Delta \omega _{1, 1} = \int_{\mathbb{T}^{d}}\frac{P\left(e^{i \boldsymbol{\theta}}\right)}{Q_{1}Q_{2}}(Q_{2}-Q_{1})dm(\boldsymbol{\theta})=0.
\label{DeltaOmega}
\end{equation}
Since $P\left(e^{i \boldsymbol{\theta}}\right) \in \mathrm{C}^{0}_{+}(\mathbb{T}^{d})$ is a strictly positive spectral density function supported on $\mathbb{T}^{d}$, and $Q_{1}, Q_{2} > 0$, the equality in \eqref{DeltaOmega} is achieved if and only if $Q_{1} = Q_{2}$. Then by Lemma \ref{Lemma43} this implies that $\Lambda_1=\Lambda_2$, establishing the injectivity of $\omega$.

Next, we prove that $\omega: \mathfrak{L}^{d}_{+} \mapsto \mathfrak{S}^{d}_{+}$ is surjective. We first note that $\omega$ is continuous and that both sets $\mathfrak{L}^{d}_{+}$ and $\mathfrak{S}^{d}_{+}$ are nonempty, convex, and open subsets of the same Euclidean space $\mathbb{T}^{(n+1)^{d}}$, and hence diffeomorphic to this space. We emphasize again that the dimension of the space is $(n+1)^{d}$ rather than $(n+1)^{2d}$, since some of the elements in $\Lambda$ and $\mathcal{T}$ are identical. 

To prove surjectivity, we will utilize Corollary 2.3 from \cite{byrnes2007interior}. According to this corollary, a continuous map \(\omega\) is surjective if and only if it is injective and proper. This means that the preimage \(\omega^{-1}(\mathfrak{K})\) must be compact for any compact subset \(\mathfrak{K}\) in \(\mathfrak{S}^{d}_{+}\). For a broader context, refer to Theorem 2.1 in \cite{byrnes2007interior}.

We now need to prove that $\omega$ is proper. To this end, we first note that $\omega^{-1}(\mathfrak{K})$ must be bounded, since if $\|\Lambda\|\to \infty$, $\omega(\Lambda)$ would tend to an all-zero matrix, which lies outside $\mathfrak{S}^{d}_{+}$ obviously. 

Consider a Cauchy sequence in $\mathfrak{K}$, which converges to a point in $\mathfrak{K}$. We need to show that the inverse image of this sequence is compact. If $\mathfrak{K}$ is empty or finite, compactness is automatic. Thus, we consider the case where the $\mathfrak{K}$ is infinite.

Since $\omega^{-1}(\mathfrak{K})$ is bounded, there must be a subsequence $(\Psi_k)$ in $\omega^{-1}(\mathfrak{K})$ that converges to a point $\Psi\in\mathfrak{L}^{d}_{+}$. It remains to show that $\Psi\in\omega^{-1}(\mathfrak{K})$, i.e., $(\Psi_{k})$ does not converge to a boundary point, which would be where $Q(e^{i\boldsymbol{\theta}})=0$. 

However, this does not occur because if it did, then $\operatorname{det} \omega(\Lambda) \rightarrow \infty$, contradicting the boundedness of $\omega^{-1}(\mathfrak{K})$. Hence, $\omega$ is proper. Therefore, the map $\omega: \mathfrak{L}^{d}_{+} \rightarrow \mathfrak{S}^{d}_{+}$ is homeomorphic.

\end{proof}

Moreover, the dual functional has the following property.
\begin{lemma}
The dual functional $\mathbb{J}_{P}(\Lambda)$ is strictly convex.
\label{Lemma45}
\end{lemma}
\begin{proof}
This is equivalent to $\delta^{2} \mathbb{J}_{P} > 0$ where
\begin{equation}
\begin{aligned}
& \delta^{2} \mathbb{J}_{P}(\Lambda; \delta \Lambda)\\
= & \int_{\mathbb{T}^{d}} \frac{P(e^{i\boldsymbol{\theta}})}{Q(e^{i\boldsymbol{\theta}}, \Lambda)^{2}}\left(K^{H}(e^{i\boldsymbol{\theta}})\delta\Lambda K(e^{i\boldsymbol{\theta}})\right)^{2}dm(\boldsymbol{\theta}).
\end{aligned}
\label{SecondDeriv}
\end{equation}

By \eqref{SecondDeriv}, we have $\delta^{2} \mathbb{J}_{P} \geqslant 0$, so it remains to show that
$$
    \delta^{2} \mathbb{J}_{P} > 0, \quad \text{for all $\delta \Lambda \neq \boldsymbol{0}$},
$$
which follows directly from Lemma \ref{Lemma43}, replacing $\Lambda$ by $\delta\Lambda$.
\end{proof}

By Lemma \ref{Lemma42}, Theorem \ref{Theorem43} and Lemma \ref{Lemma45}, we have the following theorem.

\begin{theorem}
The functional $\mathbb{J}_{P}(\Lambda)$ has a unique minimum $\Lambda^{*} \in\mathfrak{L}^{d}_{+}$. Moreover
$$
\Gamma\left(\frac{P(e^{i\boldsymbol{\theta}})}{K^{H}(e^{i\boldsymbol{\theta}})\Lambda^{*}K(e^{i\boldsymbol{\theta}})}\right)=\mathcal{T}_{d}.
$$
\label{Thm45}
\end{theorem}

By this theorem, 
$$
\Phi^{*}=\frac{P}{Q^{*}}, \quad Q^{*}=Q ( e^{i \boldsymbol{\theta}},  \Lambda^{*})
$$
belongs to $\mathrm{C}^{0}_{+}(\mathbb{T}^{d})$ and is a stationary point of $\Phi \mapsto \mathcal{L}(\Phi, \hat{\Lambda})$, which is strictly convex. Consequently
$$
\mathcal{L}(\hat{\Phi}, \hat{\Lambda}) \leqslant \mathcal{L}(\Phi, \hat{\Lambda}), \quad \text { for all } \Phi \in \mathrm{C}^{0}_{+}(\mathbb{T}^{d})
$$
or, equivalently, since $\Gamma(\hat{\Phi})=\mathcal{T}_{d}$,

$$
\mathbb{I}_{P}(P\| \hat{\Phi}) \leqslant \mathbb{KL}(P \| \Phi)
$$
for all $\Phi \in \mathrm{C}^{0}_{+}(\mathbb{T}^{d})$ satisfying the constraint $\Gamma(\Phi)=\mathcal{T}_{d}$. The above holds with equality if and only if $\Phi=\hat{\Phi}$. Finally, a complete solution is given in the following theorem for density parametrization using the multidimensional trigonometric moments.

\begin{theorem}
Let $\Gamma_{d}$ be defined by \eqref{MomentConstraintsEq} and $\Lambda$ defined by \eqref{LambdaBreve}. Given any $P \in \overline{\mathfrak{P}}_{+} \backslash \{ 0 \}$ and any $\mathcal{T}_{d} \succ 0$ with its entries $c_{k}, k = 0, \cdots, n$ given, there is a unique $\Phi \in \mathrm{C}^{0}_{+}(\mathbb{T}^{d})$ that minimizes \eqref{Entropy}
subject to $\Gamma_{d}(\Phi)=\mathcal{T}_{d}$. It has the form
\begin{equation}
\Phi^{*}=\frac{P(e^{i\boldsymbol{\theta}})}{Q(e^{i\boldsymbol{\theta}}, \Lambda^{*}_{d})},
\label{MiniForm}
\end{equation}
where $\Lambda^{*}_{d}$ is the unique solution to the problem of minimizing $\mathbb{J}_{P}$ in \eqref{LossFunc} over all $\Lambda \in \mathfrak{L}^{d}_{+}$. 
\label{Thm46}
\end{theorem}

Consequently, the dual problem provides us with an approach to compute the unique $\Phi^{*}$ that minimizes the Kullback-Leibler distance $\mathbb{KL}(P \| \Phi)$ subject to the constraint $\Gamma_{d}(\Phi)=\mathcal{T}_{d}$. 

Moreover, by Lemma \ref{Lemma42} and Theorem \ref{Theorem43}, we have the following corollary of Theorem \ref{Thm46}.

\begin{corollary}
    $\Phi^{*}\left(e^{i \boldsymbol{\theta}}\right)$ is obtained by minimizing \eqref{LossFunc} and has the form in \eqref{MiniForm}. Then $\Phi^{*}\left(e^{i \boldsymbol{\theta}}\right) > 0$ for all $\boldsymbol{\theta} \in \mathbb{T}^{d}$ if and only if $\Lambda^{*} \in \mathfrak{L}^{d}_{+}$.
\label{PhiStarEquiv}
\end{corollary}

A multidimensional parametrization using trigonometric moments has been proposed in Theorem \ref{Thm46}. Moreover, the map from the parameters of the parametrization to the trigonometric moments is proved to be homeomorphic, which also proves the wellposedness of the TMTMP of which the index set is selected as \eqref{OmegaSet}, in the Hadamard sense \cite{dontchev2006well}.

\section{Spectral estimation by the proposed parametrization}

In this section, we consider estimating the spectral density $\Phi$ by the proposed parametrization and analyzing the statistical properties. We first introduce a convex cone which is a commonly used criterion for the existence of solution to a multidimensional trigonometric moment sequence \cite{schmudgen2017moment, ringh2018multidimensional}.

Let $r = \sum_{\boldsymbol{k} \in \Omega} r_{\boldsymbol{k}} e^{-i(\boldsymbol{k}, \boldsymbol{\theta})}\in \mathbb{C}^{(n+1)d}\setminus \{0\}$, and we have $r^{*} = \sum_{\boldsymbol{k} \in \Omega} \bar{r}_{\boldsymbol{k}} e^{i(\boldsymbol{k}, \boldsymbol{\theta})}\in \mathbb{C}^{(n+1)d}\setminus \{0\}$. Denote the sequence 
\begin{equation}
c:=\left[c_{\boldsymbol{k}} \mid \boldsymbol{k} \in \Omega \right],
\label{C+}
\end{equation}
and define the open convex cone
\begin{equation}
\mathfrak{C}_{+}:=\left\{c \mid\langle c, r^{*}r\rangle>0 \right. \text{for all} \ r \in \mathbb{C}^{(n+1)d}\setminus \{0\}\}.
\end{equation}

There exists a solution to the TMTMP if and only if $c \in \mathfrak{C}_{+}$. However, as we may have noted, it is not quite easy to verify whether $c \in \mathfrak{C}_{+}$ or not, since we need the condition to be satisfied by all $r \in \mathbb{C}^{(n+1)d}\setminus \{0\}$. In the previous section, we have proposed a parametrization for the TMTMP with $\Omega$ being selected as \eqref{OmegaSet}. Given a trigonometric moment sequence $c$ satisfying $\mathcal{T}_{d} \succ 0$, there always exists a unique $\Phi^{*}$ satisfying the moment constraints exactly. To our best knowledge, it is the first parametrization which satisfies this property, compared to other treatments (see e.g. \cite{ringh2018multidimensional}). Since the trigonometric moments of the parametrization are exactly as desired, it is then feasible to perform quantitative analyses of spectral estimation using the proposed parametrization. 

We initiate this task by proposing proper statistics for the spectral estimation. We review some previous results in \cite{ringh2018multidimensional}. Denote the number of observations as $N$. For the case $d = 1$ with the scalar stationary stochastic process denoted as $\{y(t) ; t \in \mathbb{Z}\}$, it is well-known that the biased (although asymptotically unbiased) covariance estimate
$$
\hat{c}_k^{B}=\frac{1}{N} \sum_{t=0}^{N-k-1} y_t \bar{y}_{t+k},
$$
which is based on an observation record $\left\{y_t\right\}_{t=0}^{N-1}$, yields a positive definite Toeplitz matrix, which is equivalent to $\hat{c}^{B} \in \mathfrak{C}_{+}$. In fact, these estimates correspond to the ones obtained from the periodogram estimate of the spectrum (see, e.g., \cite[Sec. 2.2]{stoica1997introduction}). On the other hand, the Toeplitz matrix of the unbiased estimate
$$
\hat{c}_{k}^{U}=\frac{1}{N-k} \sum_{t=0}^{N-k-1} y_{t} \bar{y}_{t+k}
$$
is not positive definite in general.
The same holds in higher dimensions $(d>1)$ where the observation record is $\left\{y_{\boldsymbol{t}}\right\}_{\boldsymbol{t} \in \Omega}$ \cite{ringh2018multidimensional}. 

The unbiased estimate is given by
\begin{equation}
\hat{c}_{\boldsymbol{k}}^{U}=\frac{1}{\boldsymbol{N}_{\boldsymbol{k}}} \sum_{\boldsymbol{t} \in \Omega, \boldsymbol{t+k} \in \Omega} y_{\boldsymbol{t}} \bar{y}_{\boldsymbol{t}+\boldsymbol{k}},
\label{ckestimateunbiased}
\end{equation}
where we have defined $\boldsymbol{N}_{\boldsymbol{k}} := \prod_{j=1}^{d}\left(N-k_j\right)$, and the biased estimate by
\begin{equation}
\hat{c}_{\boldsymbol{k}}^{B}=\frac{1}{\boldsymbol{N}} \sum_{\boldsymbol{t} \in \Omega, , \boldsymbol{t+k} \in \Omega} y_{\boldsymbol{t}} \bar{y}_{\boldsymbol{t}+\boldsymbol{k}}
\label{ckestimate}
\end{equation}
where we have defined $\boldsymbol{N} := N^{d}$. Moreover, we define $y_{\boldsymbol{t}}=0$ for $\boldsymbol{t} \notin \Omega$. The sequence of unbiased covariance estimates does not in general belong to $\mathfrak{C}_{+}$, but the biased covariance estimates yield $\hat{c}^{B} \in \mathfrak{C}_{+}$ also in the multidimensional setting. For the spectral estimation problem, we consider the moment constraints
$$
\hat{c}_{\boldsymbol{k}}=\int_{\mathbb{T}^d} e^{ik_{1}\theta_{1}} e^{ik_{2}\theta_{2}} \cdots e^{ik_{d}\theta_{d}} \hat{\Phi}\left(e^{i \boldsymbol{\theta}}\right) d m(\boldsymbol{\theta}).
$$
where $\hat{c}_{\boldsymbol{k}}$ could be the biased estimate $\hat{c}^{B}_{\boldsymbol{k}}$ in \eqref{ckestimate} or the unbiased estimate $\hat{c}^{U}_{\boldsymbol{k}}$ in \eqref{ckestimateunbiased} (if feasible), and $\hat{\Phi}(e^{i \boldsymbol{\theta}})$ is the spectral density that we will estimate. 

In the following part of this section, we will first prove that $\mathcal{T}_{d} \succ 0$ is equivalent to $\hat{c} \in \mathfrak{C}_{+}$, of which the elements are those of $\mathcal{T}_{d}$.
\begin{proposition}
    $\mathcal{T}_{d} \succ 0$ if and only if $\hat{c} \in \mathfrak{C}_{+}$. 
\label{Prop51}
\end{proposition}

\begin{proof}
We begin by proving the sufficiency, which is straightforward. Given a $\mathcal{T}_{d} \succ 0$, by Theorem \ref{Theorem43}, there exists a unique $\Lambda \succ 0$, and then a positive $\Phi^{*}$ having the form \eqref{MiniForm}. The existence of a positive Radon measure indicates $\hat{c} \in \mathfrak{C}_{+}$ (see e.g. \cite{schmudgen2017moment}).

Next we prove the necessity. Given $\hat{c} \in \mathfrak{C}_{+}$, we have that $\Phi^{*}$ having the form \eqref{MiniForm} is everywhere positive on $\mathbb{T}^{d}$. By Corallary \ref{PhiStarEquiv}, we have that the $\Lambda^{*}_{d}$ corresponding to $\Phi^{*}$ falls within the set $\mathfrak{L}^{d}_{+}$. By Theorem \ref{Theorem43}, we obtain $\mathcal{T}_{d} \in \mathfrak{S}^{d}_{+}$, hence $\mathcal{T}_{d} \succ 0$.
\end{proof}

\begin{remark}
    In the proof of Proposition \ref{Prop51}, we have claimed that given $\hat{c} \in \mathfrak{C}_{+}$, $\Phi^{*}$ having the form \eqref{MiniForm} is everywhere positive on $\mathbb{T}^{d}$. It is not always true for a density estimate through optimization by other treatments in the literature, given the positive multidimensional moment sequence. However, by our proposed algorithm in Theorem \ref{Thm46}, the multidimensional trigonometric moments of $\Phi^{*}$ are exactly the components of $\hat{c}$. Therefore, $\Phi^{*}$ is indeed one of the members of multidimensional spectral densities which satisfy the trigonometric moments $\hat{c} \in \mathfrak{C}_{+}$. In addition, it has the least Kullback-Leibler distance from the reference spectral density $P(e^{i\boldsymbol{\theta}})$.
\end{remark}

As was mentioned in the beginning of this section, $\hat{c}^{B} \in \mathfrak{C}_{+}$ in the multidimensional setting. Then we have the corresponding matrix $\mathcal{T}_{d}$ being positive definite. However, it was usually difficult to decide whether $\hat{c}_{\boldsymbol{k}}^{U} \in \mathfrak{C}_{+}$ or not. Proposition \ref{Prop51} provides a quite straightforward criterion of evaluating the positiveness of the multidimensional trigonometric moment sequence $\hat{c}_{\boldsymbol{k}}^{U}$, with the choice of the basis functions being \eqref{OmegaSet}. When $\hat{c}_{\boldsymbol{k}}^{U} \in \mathfrak{C}_{+}$, we could then use it for the aim of spectral density estimation. 

In the following part of this section, we will then propose to analyze the statistical properties of the estimator. We first introduce the notion of Shannon-entropy maximizing spectral density. In this section, for the ease of comparing between the true density and the density estimates, we normalize the zeroth-order trigonometric moments of the true spectral density and the estimate to $1$ without loss of generality, namely,
\begin{equation}
    \int_{\mathbb{T}^d} \Phi\left(e^{i \boldsymbol{\theta}}\right) d m(\boldsymbol{\theta}) = \int_{\mathbb{T}^d} \hat{\Phi}\left(e^{i \boldsymbol{\theta}}\right) d m(\boldsymbol{\theta}) = 1.
\label{NormPhi}
\end{equation}

\begin{definition}
The Shannon-entropy maximizing spectral density $\breve{\Phi}\left(e^{i \boldsymbol{\theta}}\right)$ has the form
$$
\breve{\Phi}\left(e^{i \boldsymbol{\theta}}\right)=\exp \left(-1 - \lambda_{n+1} - \sum_{k=0}^n \lambda_k e^{i k \boldsymbol{\theta}}\right).
$$

The Lagrangian multipliers $\lambda_{0}, \cdots, \lambda_{n+1}$ are determined by maximizing the dual functional
$$
\begin{aligned} & \mathcal{L}(\Phi, \lambda)\\ = & -H(\Phi)+\sum_{k=0}^n \lambda_k\left(\int_{\mathcal{I}} e^{i k \boldsymbol{\theta}} \Phi\left(e^{i \boldsymbol{\theta}}\right) d m(\boldsymbol{\theta})-r_k\right)\\
+ & \lambda_{n+1} \Phi\left(e^{i \boldsymbol{\theta}}\right).
\end{aligned}
$$
\end{definition}

For detailed derivations of the Shannon-entropy maximization and the corresponding Shannon-entropy maximizing spectral density, please refer to our previous work \cite{wu2023quantitative}. With the proposed Shannon-entropy maximizing spectral density, we propose to derive an upper bound of difference between any two arbitrary spectral densities, which have an identical positive moment sequence $c_{\boldsymbol{k}}$.

\begin{proposition}
    Denote the set of spectral densities (of which the zeroth-order trigonometric moment is normalized to $1$), of which the multidimensional trigonometric moment sequence is $c$, as $\mathfrak{D}_{c}$. Let the total variation distance between spectral densities $\Phi$ and $\Psi$ be $V(\Phi, \Psi)$, which reads
    $$
    V(\Phi, \Psi)= \frac{1}{2}\int_{\mathbb{T}^{d}}\left|\Phi\left(e^{i \boldsymbol{\theta}}\right)-\Psi\left(e^{i \boldsymbol{\theta}}\right)\right| d m(\boldsymbol{\theta}).
    $$
    
    Given any positive multidimensional trigonometric moment sequence $c$ and two arbitrary spectral densities $\Phi, \Psi \in \mathfrak{D}_{c}$, the maximal total variation distance between $\Phi$ and $\Psi$ is upper bounded by
    \begin{equation}
    \begin{aligned}
        & \operatorname{max} V(\Phi, \Psi) = V(\breve{\Phi}, \Phi) + V(\breve{\Phi}, \Psi)\\
        \leqslant & \operatorname{min}\left\{3\left[-1+\left\{1+\frac{4}{9}(H[\breve{\Phi}]-H[\Phi])\right\}^{1 / 2}\right]^{1 / 2}\right.\\
        + & \left.3\left[-1+\left\{1+\frac{4}{9}(H[\breve{\Phi}]-H[\Psi])\right\}^{1 / 2}\right]^{1 / 2}, 1\right\},
    \end{aligned}
    \label{maxVPhiPsi}
    \end{equation}
where $\breve{\Phi}$ is the Shannon-entropy maximizing spectral density corresponding to the trigonometric moment sequence $c$.
\label{Prop63}
\end{proposition}

\begin{proof}
    We first note that the total variation distance
    $$
    \begin{aligned}
        & V(\Phi, \Psi)\\
        = & \frac{1}{2}\int_{\mathbb{T}^{d}}\left|\Phi\left(e^{i \boldsymbol{\theta}}\right)-\Psi\left(e^{i \boldsymbol{\theta}}\right)\right| d m(\boldsymbol{\theta})\\
        = & \frac{1}{2}\int_{\mathbb{T}^{d}}\left|\Phi\left(e^{i \boldsymbol{\theta}}\right) - \breve{\Phi}\left(e^{i \boldsymbol{\theta}}\right) + \breve{\Phi}\left(e^{i \boldsymbol{\theta}}\right) -\Psi\left(e^{i \boldsymbol{\theta}}\right)\right| d m(\boldsymbol{\theta})\\
        \leqslant & \frac{1}{2}\int_{\mathbb{T}^{d}}\left|\breve{\Phi}\left(e^{i \boldsymbol{\theta}}\right)-\Phi\left(e^{i \boldsymbol{\theta}}\right)\right| d m(\boldsymbol{\boldsymbol{\theta}})\\
        + & \frac{1}{2} \int_{\mathbb{T}^{d}}\left|\breve{\Phi}\left(e^{i \boldsymbol{\theta}}\right)-\Psi\left(e^{i \boldsymbol{\theta}}\right)\right| d m(\boldsymbol{\theta})\\
        = & V(\breve{\Phi}, \Phi) + V(\breve{\Phi}, \Psi).
    \end{aligned}
    $$
    The equality is achieved if and only if either of $\Phi$ and $\Psi$ is $\breve{\Phi}$, or both of them are $\breve{\Phi}$. On the other hand, by \cite{kullback1970correction, tagliani2003note}, we have
    \begin{equation}
    V(\breve{\Phi}, \Phi) \leqslant 3\left[-1+\left\{1+\frac{4}{9}(H[\breve{\Phi}]-H[\Phi])\right\}^{1 / 2}\right]^{1 / 2},
    \label{TotalVariationDis}
    \end{equation}
    and
    $$
    V(\breve{\Phi}, \Psi) \leqslant 3\left[-1+\left\{1+\frac{4}{9}(H[\breve{\Phi}]-H[\Psi])\right\}^{1 / 2}\right]^{1 / 2}.
    $$
    Moreover, we have
    $$
    \begin{aligned}
        & V(\Phi, \Psi)\\
        = & \frac{1}{2}\int_{\mathbb{T}^{d}}\left|\Phi\left(e^{i \boldsymbol{\theta}}\right)-\Psi\left(e^{i \boldsymbol{\theta}}\right)\right| d m(\boldsymbol{\theta})\\
        \leqslant & \frac{1}{2}\int_{\mathbb{T}^{d}}\Phi\left(e^{i \boldsymbol{\theta}}\right) d m(\boldsymbol{\theta}) + \frac{1}{2}\int_{\mathbb{T}^{d}}\Psi\left(e^{i \boldsymbol{\theta}}\right) d m(\boldsymbol{\theta})\\
        = & 1.
    \end{aligned}
    $$
    The third equality is because that the zeroth-order trigonometric moment of each element of the set $\mathfrak{D}_{c}$ is normalized to $1$. We have now completed the proof.
\end{proof}

Next, we will proceed to prove the consistency of the proposed estimator, considering both the biased trigonometric moments $\hat{c}^{B}$ and the unbiased trigonometric moments $\hat{c}^{U}$. 

\begin{theorem}
    With the number of observations, namely $N$, and the prescribed highest order of the trigonometric moments, namely $n$, approaching infinity, $V(\hat{\Phi}, \Phi)$ converges to $0$, i.e.,  the density estimate $\hat{\Phi}\left(e^{i \boldsymbol{\theta}}\right)$ is equal to the true spectral density $\Phi\left(e^{i \boldsymbol{\theta}}\right)$, with the statistics being either $\hat{c}^{B}$ or $\hat{c}^{U}$. Furthermore, with $N \rightarrow +\infty$, the upper bound of $V(\hat{\Phi}, \Phi)$, which has the form in \eqref{maxVPhiPsi}, decreases monotonically with the increase of $n$.
\label{Thrm44}
\end{theorem}

\begin{proof}
Considering the biased trigonometric moments $\hat{c}^{B}$, we have that with the number of the observations $N \rightarrow +\infty$, the biased covariance estimates $\hat{c}^{B} \rightarrow c$. The same holds for the unbiased covariance estimates $\hat{c}^{U}$. Therefore, in the following part of this proof, we consider $V(\hat{\Phi}, \Phi)$ where $\hat{\Phi}$ is obtained by the true trigonometric moment sequence $c$.

With $n \rightarrow +\infty$, $\hat{\Phi}$ and $\Phi$ have the identical full trigonometric moment sequence, which leads to the fact that 
$$
\hat{\Phi} = \Phi, \text{a.e.}.
$$
In the problem setting of this paper we have confined $\hat{\Phi}$ and $\Phi$ to fall within the set $\mathrm{C}^{0}_{+}(\mathbb{T}^{d})$, therefore we have $\hat{\Phi} = \Phi$ (see e.g. \cite[Theorem 11.3]{schmudgen2017moment} for more details), which proves the consistency. 

Next we prove that the upper bound of $V(\hat{\Phi}, \Phi)$, with $N \rightarrow +\infty$, decreases monotonically with the increase of $n$. With a slight abuse of use of the symbols, we denote $\hat{\Phi}_{[m]}$ as the spectral density estimate by $c_{[m]}$, which is the true trigonometric moment sequence of moments up to order $m$. Given arbitrary $m_{1}, m_{2} \in \mathbb{N}, m_{1} < m_{2}$, we have
\begin{equation}
\mathfrak{D}_{c_{[m_{2}]}} \subseteq \mathfrak{D}_{c_{[m_{1}]}}.
\label{Dsubset}
\end{equation}
It is obvious since for a longer trigonometric moment sequence, higher order moment specifications need to be satisfied. Moreover, we note that
\begin{equation}
H\left[\breve{\Phi}_{[m_{2}]}\right] \leqslant H\left[\breve{\Phi}_{[m_{1}]}\right].
\label{Hm1m2}
\end{equation}
It is also quite obvious since $\breve{\Phi}_{[m_{2}]}$ is the one with largest Shannon-entropy in $\mathfrak{D}_{c_{[m_{2}]}}$, and $\breve{\Phi}_{[m_{1}]}$ is the one with largest Shannon-entropy in $\mathfrak{D}_{c_{[m_{1}]}}$. Therefore, \eqref{Dsubset} naturally leads to \eqref{Hm1m2}, which completes the proof. We also note that with $n \rightarrow +\infty$, $\Phi$, $\breve{\Phi}$, $\hat{\Phi}$ have the identical full trigonometric moment sequence, hence $\operatorname{max} V(\Phi, \hat{\Phi}) \rightarrow 0$, which also proves the consistency of the estimator in the sense of the total variation distance.
\end{proof}

Theorem \ref{Thrm44} proves the consistency of the proposed spectral density estimator. Next, we prove the (asymptotic) unbiasedness of the proposed estimator, given that $n \rightarrow +\infty$. 

\begin{proposition}
    Given that $n \rightarrow +\infty$, the spectral density estimate, obtained by Theorem \ref{Thm46}, is unbiased, if the sample covariance sequence used is $\hat{c}^{U}$. It is an asymptotically unbiased estimate if the sample covariance sequence $\hat{c}^{B}$ is utilized.
\label{Prop55}
\end{proposition}

\begin{proof}
   Denote $\Omega_{\infty}$ as the $\Omega$ with $n \rightarrow +\infty$. We note that the multidimensional spectral density $\Phi$ of a purely nondeterministic, second-order, stationary random process $\{y_{\boldsymbol{t}}\}$ with zero mean is given by the Fourier expansion
    \begin{equation}
        \Phi\left(e^{i \boldsymbol{\theta}}\right) = \sum_{k \in \Omega_{\infty}} c_{\boldsymbol{k}} e^{i(\boldsymbol{k}, \boldsymbol{\theta})}.
    \label{PhiInfExpansion}
    \end{equation}
    With the unbiased truncated sample covariance sequence $\hat{c}^{U}$, the spectral estimate by Theorem \ref{Thm46} could be written as
    $$
        \hat{\Phi}\left(e^{i \boldsymbol{\theta}}\right) = \sum_{\boldsymbol{k} \in \Omega} \hat{c}^{U}_{\boldsymbol{k}} e^{i(\boldsymbol{k}, \boldsymbol{\theta})} + \sum_{\boldsymbol{k} \in \Omega_{\infty} \setminus \Omega} \tilde{c}^{U}_{\boldsymbol{k}} e^{i(\boldsymbol{k}, \boldsymbol{\theta})},
    $$
    where
    $$
        \tilde{c}^{U}_{\boldsymbol{k}} = \int_{\mathbb{T}^d} e^{i(\boldsymbol{k}, \boldsymbol{\theta})} \hat{\Phi}\left(e^{i \boldsymbol{\theta}}\right) d m
    $$
    with $\hat{\Phi}$ being estimated by Theorem \ref{Thm46}. We note that with $n \rightarrow +\infty$, the estimates of the frequency components of the density span the whole domain of $\boldsymbol{k}$, namely the $d$-dimensional real space, which leads to
    \begin{equation}
        \operatorname{lim}_{n \rightarrow +\infty} \hat{\Phi}(e^{i \boldsymbol{\theta}}) = \sum_{\boldsymbol{k} \in \Omega_{\infty}} \hat{c}^{U}_{\boldsymbol{k}} e^{i(\boldsymbol{k}, \boldsymbol{\theta})}.
    \label{PhiEstInf}
    \end{equation}
    Hence to prove the unbiasedness of the spectral estimate $\hat{\Phi}$ is equivalent to prove that each frequency component of the spectral estimate, namely $\hat{c}^{U}_{\boldsymbol{k}}$ is unbiased. This condition is satisfied since $\hat{c}^{U}_{\boldsymbol{k}}$ is an unbiased estimate of $c_{\boldsymbol{k}}$. Similarly, we could prove that the spectral density estimate $\hat{\Phi}$ by the biased sample covariance sequence is asymptotically unbiased since each frequency component $\hat{c}^{B}_{\boldsymbol{k}}$ of the spectral estimate is an asymptotically unbiased estimate of $c_{\boldsymbol{k}}$. Now the proof is complete.
\end{proof}

We would also like to investigate the convergence rate of the spectral estimate $\hat{\Phi}(e^{i\boldsymbol{\theta}})$. We note that the sample trigonometric moments $\hat{c}^{U}_{\boldsymbol{k}}$ and $\hat{c}^{B}_{\boldsymbol{k}}$ are an unbiased estimate and an asymptotically unbiased estimate of $c_{\boldsymbol{k}}$ respectively, for which we assume the variance as $\operatorname{var}(\hat{c}^{U}_{\boldsymbol{k}}) = \operatorname{var}(\hat{c}^{B}_{\boldsymbol{k}}) =\sigma_{\boldsymbol{k}}^2$ (It is quite straightforward to prove the first equality and we omit the proof here). The following results hold for both $\hat{c}^{U}_{\boldsymbol{k}}$ and $\hat{c}^{B}_{\boldsymbol{k}}$, hence we adopt the notation $\hat{c}_{\boldsymbol{k}}$ for the following derivations. For a large $\boldsymbol{N}_{\boldsymbol{k}}$ we approximately have $\sqrt{\boldsymbol{N}_{\boldsymbol{k}}}(\hat{c}_{\boldsymbol{k}}-c_{\boldsymbol{k}}) \sim \mathcal{N}\left(0, \sigma_{\boldsymbol{k}}^2\right)$ by the central limit theorem (we emphasize that $\boldsymbol{N}_{\boldsymbol{k}}$ is valid for both $\hat{c}_{\boldsymbol{k}}^{U}$ and $\hat{c}_{\boldsymbol{k}}^{B}$ since they are both sums of $\boldsymbol{N}_{\boldsymbol{k}}$ terms by \eqref{ckestimateunbiased} and \eqref{ckestimate}). Then for any $\epsilon_{\boldsymbol{k}} > 0$, we have that there exist unique constants $\alpha_{\boldsymbol{k}} > 0, \boldsymbol{k} \in \Omega$ satisfying
\begin{equation}
\mathbb{P}\left(\left|\hat{c}_{\boldsymbol{k}}-c_{\boldsymbol{k}}\right| \geqslant \alpha_{\boldsymbol{k}} \sigma_{\boldsymbol{k}} \cdot \boldsymbol{N}_{\boldsymbol{k}}^{-1 / 2}\right) \leqslant \epsilon_{\boldsymbol{k}},
\label{OpDef}
\end{equation}
by which we can conclude that 
\begin{equation}
\hat{c}_{\boldsymbol{k}}-c_{\boldsymbol{k}}=\mathrm{O}_\mathbb{P}\left(\boldsymbol{N}_{\boldsymbol{k}}^{-1 / 2}\right)
\label{OpMoment}
\end{equation}
for $\boldsymbol{k} \in \Omega$, where the definition of $\mathrm{O}_\mathbb{P}(\cdot)$ can be found in \cite{dasgupta2008asymptotic}. On the other hand, for any $v>1 / 2$ we have $\boldsymbol{N}_{\boldsymbol{k}}^{-v} / \boldsymbol{N}_{\boldsymbol{k}}^{-1 / 2} \rightarrow 0$ as $\boldsymbol{N}_{\boldsymbol{k}} \rightarrow \infty$. Hence, for any constant $\alpha_{\boldsymbol{k}}>0$,
$$
\begin{aligned}
& \mathbb{P} \left(\left|\hat{c}_{\boldsymbol{k}}-c_{\boldsymbol{k}}\right| \geqslant \alpha_{\boldsymbol{k}} \sigma_{\boldsymbol{k}} \cdot \boldsymbol{N}_{\boldsymbol{k}}^{-v}\right) \\
= & \mathbb{P} \left(\left|\hat{c}_{\boldsymbol{k}}-c_{\boldsymbol{k}}\right| \geqslant \left(\alpha_{\boldsymbol{k}} \sigma_{\boldsymbol{k}} \cdot \boldsymbol{N}_{\boldsymbol{k}}^{-1 / 2}\right) \cdot \frac{\boldsymbol{N}_{\boldsymbol{k}}^{-v}}{\boldsymbol{N}_{\boldsymbol{k}}^{-1 / 2}}\right)\\
\rightarrow & 1 \text { as } \boldsymbol{N}_{\boldsymbol{k}} \rightarrow \infty,
\end{aligned}
$$
for $\boldsymbol{k} \in \Omega$, which reveals the fact that with $\boldsymbol{N}_{\boldsymbol{k}} \rightarrow +\infty$, \eqref{OpDef} cannot be achieved. Therefore $\boldsymbol{N}_{\boldsymbol{k}}^{-1 / 2}$ is the convergence rate of $\hat{c}_{\boldsymbol{k}}$.

By \eqref{PhiInfExpansion}, we further have
\begin{equation}
\begin{aligned}
& E_{\Phi}\left[ \hat{\Phi} - \Phi \right]\\
= & \int_{\mathbb{T}^{d}} \left( \hat{\Phi}(e^{i \boldsymbol{\theta}}) - \Phi(e^{i \boldsymbol{\theta}}) \right) \cdot \Phi(e^{i \boldsymbol{\theta}}) d m(\boldsymbol{\theta})\\
= & \sum_{k \in \Omega_{\infty}} c_{\boldsymbol{k}}\int_{\mathbb{R}} e^{i(\boldsymbol{k}, \boldsymbol{\theta})} \cdot \left( \hat{\Phi}(e^{i \boldsymbol{\theta}}) - \Phi(e^{i \boldsymbol{\theta}}) \right) d m(\boldsymbol{\theta})\\
= & \sum_{k \in \Omega_{\infty}} c_{\boldsymbol{k}} \left( \hat{c}_{\boldsymbol{k}} - c_{\boldsymbol{k}} \right)
\end{aligned}
\label{EpSub}
\end{equation}

Since $\hat{\Phi} \in \mathrm{C}^{0}_{+}(\mathbb{T}^{d})$ and is finite, we assume that there exists a $\beta \in \mathbb{R}_{+}$ such that $\hat{\Phi} \leqslant \beta$ on $\mathbb{T}^{d}$. By \eqref{NormPhi} we further have
$$
\begin{aligned}
& E_{\Phi}\left[ \hat{\Phi} - \Phi \right]\\
= & \int_{\mathbb{T}^{d}} \left( \hat{\Phi} - \Phi \right) \cdot \Phi d m(\boldsymbol{\theta})\\
\leqslant & \int_{\mathbb{T}^{d}} \hat{\Phi} \cdot \Phi d m(\boldsymbol{\theta})\\
\leqslant & \int_{\mathbb{T}^{d}} \beta \Phi d m(\boldsymbol{\theta}) = \beta.
\end{aligned}
$$
Hence $E_{\Phi}\left[ \hat{\Phi} - \Phi \right]$ is finite.

We note that each term on the RHS of the last equation of \eqref{EpSub} satisfies 
$$
c_{\boldsymbol{k}} \left( \hat{c}_{\boldsymbol{k}} - c_{\boldsymbol{k}} \right) = \mathrm{O}_\mathbb{P}\left(\boldsymbol{N}_{\boldsymbol{k}}^{-1 / 2}\right)
$$
for $\boldsymbol{k} \in \Omega_{\infty} \setminus \Omega$, and $E_{\Phi}\left[ \hat{\Phi} - \Phi \right]$ is finite. We can then conclude that $E_{\Phi}\left[ \hat{\Phi} - \Phi \right] = \mathrm{O}_\mathbb{P}\left(\boldsymbol{N}_{\boldsymbol{k}}^{-1 / 2}\right)$.

Therefore, we prove the uniform convergence and the convergence rate of the proposed density estimator to be $\boldsymbol{N}_{\boldsymbol{k}}^{-1/2}$, which is the optimal (fastest possible) convergence rate of parametric estimation \cite{silverman2018density}. On the other hand, in nonparametric density estimation, the optimal convergence rate is $\boldsymbol{N}_{\boldsymbol{k}}^{-s/(2s+1)}$ only if the underlying function is continuously differentiable up to order $s$ \cite{he1994convergence}. The convergence rate highlights the superiority of the proposed spectral estimator compared to other nonparametric estimators. Moreover, since it doesn't require a prior assumption of the spectral density, the proposed density parametrization by multidimensional trigonometric moments provides much more flexibility, compared to conventional parametric methods. 

Moreover, we could further prove the efficiency of the proposed spectral density estimate. For the proof of efficiency, we would first need to prove the unbiasedness of the estimator, which is given in Theorem \ref{Prop55}. Again by the central limit theorem, we note that given a large $\boldsymbol{N}_{\boldsymbol{k}}$, the estimate of the covariance lag $\hat{c}_{\boldsymbol{k}} \sim \mathcal{N}(c_{\boldsymbol{k}}, \frac{\sigma^{2}_{\boldsymbol{k}}}{\boldsymbol{N}_{\boldsymbol{k}}})$, $\boldsymbol{k} \in \Omega$ by the central limit theorem. To prove the efficiency, it suffices to prove that the variance of $\hat{c}_{\boldsymbol{k}}$ is equal to the Cramér-Rao Lower Bound (CRLB), which provides a lower bound on the variance of any unbiased estimator of a parameter. Now we prove the following proposition.

\begin{proposition}
    Assume that the stationary stochastic process $\{y_{\boldsymbol{t}}\}$ satisfies mild mixing conditions (or is a Gaussian process) such that its sample covariances satisfy the multivariate central limit theorem. Then, the proposed spectral density estimate $\hat{\Phi}(e^{i\boldsymbol{\theta}})$ is an asymptotically efficient estimate of the true spectral density function $\Phi(e^{i\boldsymbol{\theta}})$ as the sample size $N \to \infty$ and the basis dimension $n \to \infty$.
\end{proposition}

\begin{proof}
    Instead of assuming the individual lag products $c_{\boldsymbol{k, t}} := y_{\boldsymbol{t}} \bar{y}_{\boldsymbol{t}+\boldsymbol{k}}$ are independent and normally distributed, which is generally invalid since products of dependent random variables are non-Gaussian, we establish the asymptotic efficiency by invoking the asymptotic normality of the aggregate sample covariance vector and the smooth mapping property of our optimization framework.

    Let $\hat{\boldsymbol{c}}_n = [\hat{c}_{\boldsymbol{k}}]_{\boldsymbol{k} \in \Omega}$ and $\boldsymbol{c}_n = [c_{\boldsymbol{k}}]_{\boldsymbol{k} \in \Omega}$ denote the vector of sample covariances and true trigonometric moments, respectively. Under mild regularity and mixing conditions for stationary random fields (see, e.g., \cite{brockwell2009time}), the sample covariance vector $\hat{\boldsymbol{c}}_n$ satisfies the multivariate Central Limit Theorem (CLT) as the total number of available observation samples $N \to \infty$:
    \begin{equation}
        \sqrt{N} (\hat{\boldsymbol{c}}_n - \boldsymbol{c}_n) \xrightarrow{d} \mathcal{N}(\mathbf{0}, \boldsymbol{\Sigma}_n),
    \end{equation}
    where $\xrightarrow{d}$ denotes convergence in distribution, and $\boldsymbol{\Sigma}_n$ is the asymptotic covariance matrix determined by the higher-order moments of the process $y_{\boldsymbol{t}}$. When $\{y_{\boldsymbol{t}}\}$ is a Gaussian process, the sample covariance vector $\hat{\boldsymbol{c}}_n$ is known to be an asymptotically efficient estimator, meaning its asymptotic covariance matrix $\frac{1}{N}\boldsymbol{\Sigma}_n$ coincides with the minimum variance bound defined by the CRLB for the covariance parameters.

    Recall that from the variational analysis in Section~\ref{sec3}, the proposed convex optimization scheme establishes a diffeomorphism between the feasible moment sequence $\boldsymbol{c}_n$ and the spectral parameter vector $\Lambda$. Consequently, the spectral density estimate can be viewed as a continuously differentiable (smooth) functional of the sample moments, denoted by $\hat{\Phi}(e^{i\boldsymbol{\theta}}) = \mathcal{G}(\hat{\boldsymbol{c}}_n)$. 

    By applying the Multivariate Delta Method \cite{van2000asymptotic}, the asymptotic distribution of the spectral density estimate inherits the normality and efficiency of the underlying moment estimator. Specifically, for any fixed frequency vector $\boldsymbol{\theta}$ within the multidimensional torus $\mathbb{T}^d$, the pointwise asymptotic distribution satisfies
\begin{equation}
    \sqrt{N} \left(\hat{\Phi}(e^{i\boldsymbol{\theta}}) - \Phi(e^{i\boldsymbol{\theta}})\right) \xrightarrow{d} \mathcal{N}\left(0, \nabla \mathcal{G}_{\boldsymbol{\theta}}(\boldsymbol{c}_n)^T \boldsymbol{\Sigma}_n \nabla \mathcal{G}_{\boldsymbol{\theta}}(\boldsymbol{c}_n)\right),
\end{equation}
where $\mathcal{G}_{\boldsymbol{\theta}}(\cdot)$ denotes the smooth functional evaluating the spectral density at $\boldsymbol{\theta}$ from the moment sequence, and $\nabla \mathcal{G}_{\boldsymbol{\theta}}(\boldsymbol{c}_n)$ represents its gradient evaluated at the true moments $\boldsymbol{c}_n$. 
Since $\hat{\boldsymbol{c}}_n$ asymptotically achieves the CRLB and the mapping $\mathcal{G}_{\boldsymbol{\theta}}$ is a bijective diffeomorphism within the interior of the feasible domain, the transformed estimator $\hat{\Phi}(e^{i\boldsymbol{\theta}})$ also asymptotically achieves the corresponding information theoretic lower bound for the spectral density function pointwise for a fixed basis dimension $n$ \cite{lehmann1998theory}.

Finally, as the truncation order and basis dimension $n \rightarrow +\infty$, the chosen basis functions form a complete set in the continuous function space. According to \eqref{PhiInfExpansion}, the approximation error uniformly vanishes. Combining the statistical efficiency under $N \to \infty$ with the structural completeness under $n \to \infty$, we conclude via \eqref{PhiEstInf} that $\hat{\Phi}(e^{i\boldsymbol{\theta}})$ is an asymptotically efficient estimate of the true multidimensional spectral density function.
\end{proof}

In this section, we have proposed an estimation scheme of the spectral density, together with a sufficient and necessary condition for the estimate to be positive. Statistical properties of the estimate, including the consistency, (asymptotic) unbiasedness, the convergence rate, and the asymptotic efficiency under a mild assumption are proved. Based on these observations, we propose a comprehensive algorithm to address the multidimensional rational covariance extension problem. However, there still remains an interesting problem which draws our interest. By Theorem \ref{Thm46}, we have that given any $\theta \in \mathcal{P}$, there exists a unique $\Phi^{*}$, corresponding to a unique $\Lambda^{*}_{d}$, minimizing the KL distance between a reference distribution and $\Phi^{*}$. In some scenarios, we would like to know the modes (peaks) of the density estimate. The modes of a spectral density function represent the dominant frequencies present in the wide-sense stationary process. These peaks provide crucial information about the periodic components and could help identify underlying patterns or oscillatory behavior in the data. However, in the literature of unidimensional or multidimensional RCEP, there has not yet been a result on the number of the modes of the spectral estimate and the prior spectral density $P(e^{i \boldsymbol{\theta}})$. In the following part of this section, we will analyze the modes of a spectral estimate, given a specific choice of $P(e^{i \boldsymbol{\theta}})$.

\begin{proposition}
Define 
$$
\Omega_{+} := \left\{ \left(k_{1}, \cdots k_{d} \right) \mid 0 \leqslant k_{i} \leqslant n, n \in \mathbb{N}, i = 1, \cdots, d \right\}.
$$
Let the spectral estimate $\hat{\Phi}(e^{i \boldsymbol{\theta}}) = \frac{P(e^{i \boldsymbol{\theta}})}{Q(e^{i \boldsymbol{\theta}})}$ be obtained by Theorem \ref{Thm46}. Then we have that $\hat{\Phi}(e^{i \boldsymbol{\theta}})$ has up to $(n(n-1))^{d}$ modes.
\label{Prop56}
\end{proposition}

\begin{proof}
First we write $\hat{\Phi}(e^{i \boldsymbol{\theta}})$ 
\begin{equation}
\hat{\Phi}(e^{i \boldsymbol{\theta}}) = w(e^{i \boldsymbol{\theta}})w(e^{-i \boldsymbol{\theta}}) = |w(e^{i \boldsymbol{\theta}})|^{2},
\label{Wsquare}
\end{equation}
following \cite{byrnes2002identifiability} where
$$
w(e^{i \boldsymbol{\theta}}) = \frac{a(e^{i\boldsymbol{\theta}})}{b(e^{i\boldsymbol{\theta}})} = \frac{\sum_{\boldsymbol{k} \in \Omega_{+}} a_{\boldsymbol{k}}e^{i \boldsymbol{k} \boldsymbol{\theta}}}{\sum_{\boldsymbol{k} \in \Omega_{+}} b_{\boldsymbol{k}}e^{i \boldsymbol{k} \boldsymbol{\theta}}}.
$$
We note that $a(e^{i \boldsymbol{\theta}}), b(e^{i \boldsymbol{\theta}})$ are both trigonometric polynomials of order $n$.

Then we have that the first order derivative of $w(e^{i \boldsymbol{\theta}})$ reads
\begin{equation}
\begin{aligned}
\frac{\partial w(e^{i \boldsymbol{\theta}})}{\partial \boldsymbol{\theta}}= & \frac{\frac{\partial a(e^{i \boldsymbol{\theta}})}{\partial e^{i \boldsymbol{\theta}}} \frac{\partial e^{i \boldsymbol{\theta}}}{\boldsymbol{\partial \theta}}b(e^{i \boldsymbol{\theta}}) - \frac{\partial b(e^{i \boldsymbol{\theta}})}{\partial e^{i \boldsymbol{\theta}}} \frac{\partial e^{i \boldsymbol{\theta}}} {\boldsymbol{\partial \theta}}a(e^{i \boldsymbol{\theta}})}{b^2(e^{i \boldsymbol{\theta}})}\\
= & \frac{ i \left(\frac{\partial a(e^{i \boldsymbol{\theta}})}{\partial e^{i \boldsymbol{\theta}}} b(e^{i \boldsymbol{\theta}}) - \frac{\partial b(e^{i \boldsymbol{\theta}})}{\partial e^{i \boldsymbol{\theta}}} a(e^{i \boldsymbol{\theta}}) \right)}{b^2(e^{i \boldsymbol{\theta}})}.
\end{aligned}
\label{PartialPhiTheta}
\end{equation}

Since $a(e^{i \boldsymbol{\theta}})$ and $b(e^{i \boldsymbol{\theta}})$ are both trigonometric polynomials of order $n$, we have that $\frac{\partial a(e^{i \boldsymbol{\theta}})}{\partial e^{i \boldsymbol{\theta}}}$ and $\frac{\partial b(e^{i \boldsymbol{\theta}})}{\partial e^{i \boldsymbol{\theta}}}$ are both trigonometric polynomials of order $n-1$. Since the denominator, namely $b^2(e^{i \boldsymbol{\theta}})$ is positive, and the numerator is a trigonometric polynomial of order $n(n-1)$, we then have that $\frac{\partial w(e^{i \boldsymbol{\theta}})}{\partial \boldsymbol{\theta}}$ has at most $(n(n-1))^{d}$ zeros. Since $n \in \mathbb{N}$, $(n(n-1))^{d}$ is an even integer. Therefore, $w(e^{i \boldsymbol{\theta}})$ has at most $\frac{(n(n-1))^{d}}{2}$ modes (peaks). Then by \eqref{Wsquare}, we have that $\hat{\Phi}(e^{i \boldsymbol{\theta}})$ has at most $(n(n-1))^{d}$ modes (peaks).
\end{proof}

By Proposition \ref{Prop56}, we could then obtain the minimal $n$, namely the highest order of the trigonometric moments, given the number of modes (peaks) that we would like to have in the spectral estimate.

We then conclude the contributions that we have made to the TMTMP and the RCEP in this section. In the literature, we are always confronted with two types of multidimensional trigonometric moment estimates, namely the biased trigonometric moment estimates $\hat{c}^{B}$ and the unbiased ones $\hat{c}^{U}$. For the unidimensional RCEP, the existence of solution is ensured by the positive definiteness of the Toeplitz matrix formed by the trigonometric moments. In general, the biased trigonometric moment estimates $\hat{c}^B$ satisfy this condition, while the unbiased counterparts $\hat{c}^U$ don't necessarily do so. However, we could easily check whether a solution exists given specific $\hat{c}^{U}$ by checking the positive definiteness of the Toeplitz matrix formed by the trigonometric moments $\hat{c}^{U}$. But for the multidimensional RCEP, there has not been such an explicit condition to guarantee the existence of solution given a specific $\hat{c}^{U}$. In \cite{ringh2018multidimensional}, an estimation approach by approximate covariance matching is proposed, which seeks to estimate the spectral density by matching $\hat{c}^{U}$ to a sequence of trigonometric moment estimates in $\mathfrak{C}_{+}$, which is close to $\hat{c}^{U}$ in some norm. The simulation results show the superiority of the proposed approximate moment matching scheme. Despite the numerical performance, this treatment somehow lacks evidence in statistics, since the sequence of trigonometric moment estimates in $\mathfrak{C}_{+}$, obtained by considering some norm as a term of the cost function, doesn't necessarily reveal itself as proper statistics of the true spectral density. In this paper, on the contrary to approximating $\hat{c}^{U}$ to a feasible one in $\mathfrak{C}_{+}$ by closeness in some norm, we propose an explicit condition for $\hat{c}^{U} \in \mathfrak{C}_{+}$ in Proposition \ref{Prop51}. The positive definiteness of the matrix $\mathcal{T}_{d}$, of which the entries are those of $\hat{c}^{U}$, is now the sufficient and necessary condition of $\hat{c}^{U} \in \mathfrak{C}_{+}$. Hence the estimation of the spectral density could be performed as follows. If $\hat{c}^{U} \in \mathfrak{C}_{+}$, the corresponding $\mathcal{T}_{d} \succ 0$ and we use $\hat{c}^{U}$ for the estimation task by Theorem \ref{Thm46}; otherwise, we use $\hat{c}^{B}$ for estimation instead. A detailed algorithm is given in Theorem 1.

\begin{algorithm}[t]
	\caption{A multidimensional rational covariance extension algorithm}
	\label{alg:algorithm1}
	\KwIn{Observation record $\left\{y_{\boldsymbol{t}}\right\}_{\boldsymbol{t} \in \Omega}$ of the WSS stochastic process corresponding to $\Phi(e^{i\theta})$}
	\KwOut{Spectral density function estimate $\hat{\Phi}(e^{i\boldsymbol{\theta}})$.}  
	\BlankLine
	Step 1: Calculate $\hat{c}^{B}$ and $\hat{c}^{U}$ by \eqref{ckestimateunbiased} and \eqref{ckestimate} correspondingly;\\
	Step 2: Examine whether $\hat{c}^{U} \in \mathfrak{C}_{+}$ by Proposition \ref{Prop51}. If $\hat{c}^{U} \in \mathfrak{C}_{+}$, we take $\hat{c} = \hat{c}^{U}$; else, we take $\hat{c} = \hat{c}^{B}$ instead.\\
        Step 3: Obtain the spectral density estimate $\hat{\Phi}(e^{i\boldsymbol{\theta}})$ by Theorem \ref{Thm46} using the sequence of multidimensional trigonometric moments $\hat{c}$.
\end{algorithm}

By doing this, the proposed spectral estimate satisfies the statistical properties proposed in this section, including the consistency, (asymptotic) unbiasedness, and the convergence rate. Moreover, with $N \rightarrow +\infty$, the upper bound of $V(\hat{\Phi}, \Phi)$ decreases monotonically with the increase of $n$. It is the first time that the statistical properties of the spectral estimate are comprehensively analyzed, and the choice of $n$ is proved to be the larger the better, in the sense of the total variation distance. 

\section{A numerical example}

In the previous sections, we provided a detailed treatment of the TMTMP from a system and signal processing perspective, comprehensively analyzing the statistical properties of the proposed multidimensional spectral density estimator. In this section, we present numerical simulations to empirically validate the efficacy of the proposed approach.

\begin{figure} [htbp]
    \centering
    \subfigure[]{
    \includegraphics[scale=0.48]{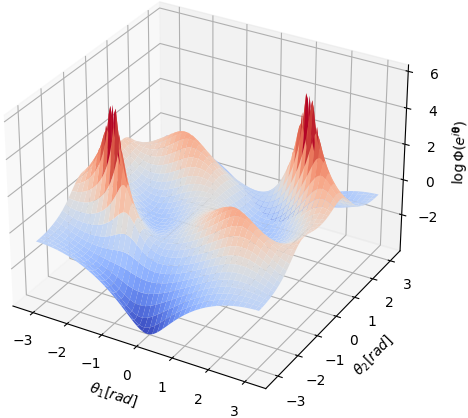}
    }
    \subfigure[]{
    \includegraphics[scale=0.48]{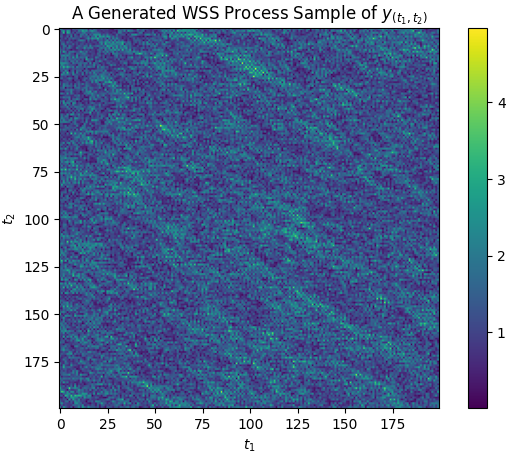}
    }
    \caption{(a) Log-plot of the original spectrum; (b) The norm of a generated wide-sense stationary process sample of $y_{(t_{1}, t_{2})}$ with $t_{1} \leqslant 200$, $t_{2} \leqslant 200$.}
    \label{Fig1}
\end{figure}

We consider the example presented in Subsection 9.1 of \cite{ringh2018multidimensional}. Let $y_{(t_1, t_2)}$ be the steady-state output of a two-dimensional recursive filter driven by a white noise input $u_{(t_1, t_2)}$. The transfer function of the recursive filter is defined as
$$
\frac{b(e^{i \theta_1}, e^{i \theta_2})}{a(e^{i \theta_1}, e^{i \theta_2})}=\frac{\sum_{\boldsymbol{k} \in \Omega_{+}} b_{\boldsymbol{k}} e^{-i(\boldsymbol{k}, \boldsymbol{\theta})}}{\sum_{\boldsymbol{k} \in \Omega_{+}} a_{\boldsymbol{k}} e^{-i(\boldsymbol{k}, \boldsymbol{\theta})}}
$$
where $\Omega_{+}=\{(k_1, k_2) \in \mathbb{Z}^2 \mid 0 \leqslant k_1 \leqslant 2, 0 \leqslant k_2 \leqslant 2\}$. Denoting the entry of a matrix $A$ in row $k_1$ and column $k_2$ as $A_{k_1, k_2}$, the coefficients are given by $b_{(k_1, k_2)}=B_{k_1+1, k_2+1}$ and $a_{(k_1, k_2)}=A_{k_1+1, k_2+1}$, with
$$
B=\left[\begin{array}{ccc}
0.9 & -0.2 & 0.05 \\
0.2 & 0.3 & 0.05 \\
-0.05 & -0.05 & 0.1
\end{array}\right], \\
$$
and
$$
A=\left[\begin{array}{crr}
1 & 0.1 & 0.1 \\
-0.2 & 0.2 & -0.1 \\
0.4 & -0.1 & -0.2
\end{array}\right] .
$$
The true spectral density $\Phi$ of $y_{(t_1, t_2)}$, depicted in Figure \ref{Fig1}(a) and structurally similar to the model analyzed in \cite{ringh2016multidimensional}, is formulated as
$$
\Phi(e^{i \theta_1}, e^{i \theta_2})=\frac{P(e^{i \theta_1}, e^{i \theta_2})}{Q(e^{i \theta_1}, e^{i \theta_2})}=\left|\frac{b(e^{i \theta_1}, e^{i \theta_2})}{a(e^{i \theta_1}, e^{i \theta_2})}\right|^2.
$$
Consequently, the index set $\Omega$ corresponding to the coefficients of the trigonometric polynomials $P$ and $Q$ is determined by the Minkowski set difference $\Omega=\Omega_{+}-\Omega_{+}=\{(k_1, k_2) \in \mathbb{Z}^2 \mid | k_1|\leqslant 2,|k_2| \leqslant 2\}$. 

Formally, within the rational covariance extension framework, the estimation objective is to identify the optimal denominator polynomial parameters $\Lambda$ given a prior numerator polynomial $P(e^{i \boldsymbol{\theta}})$. The parameterized spectral estimate is expressed as
$$
\Phi(e^{i \boldsymbol{\theta}} ; \Lambda) = \frac{P(e^{i \boldsymbol{\theta}})}{Q(e^{i \boldsymbol{\theta}}, \Lambda)} = \frac{P(e^{i \boldsymbol{\theta}})}{\sum_{\boldsymbol{k} \in \Omega} \lambda_{\boldsymbol{k}} e^{-i(\boldsymbol{k}, \boldsymbol{\theta})}}.
$$
To strictly evaluate the convexity and robustness of the proposed optimization scheme for the autoregressive (AR) parameters $\Lambda$, we fix the moving average (MA) part by prescribing $P(e^{i \boldsymbol{\theta}})$ to identically match the true numerator. This decoupling is mathematically non-trivial; as established in \cite{byrnes2002identifiability}, the mapping from the numerator coefficients to the cepstral coefficients of the spectral estimate—given a fixed denominator—constitutes a diffeomorphism, which uniquely dictates the zeroes and structural smoothness of the spectrum.

To empirically validate the estimator, sample realizations of the 2D wide-sense stationary (WSS) process are generated. This is achieved by passing a 2D standard Gaussian white noise process through the true rational filter via the inverse Fourier transform, yielding the time-domain realization $y_{(t_1, t_2)}$. A sample trajectory is illustrated in Figure \ref{Fig1}(b). The sample statistics of the multidimensional trigonometric moments are subsequently computed via \eqref{ckestimateunbiased} and \eqref{ckestimate}.

\begin{figure}[htbp]
    \centering
    \includegraphics[width=0.48\textwidth]{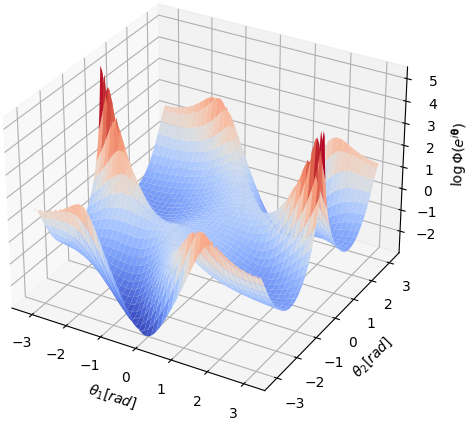}
    \includegraphics[width=0.48\textwidth]{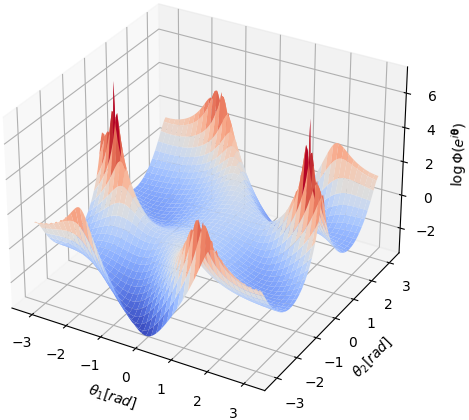}\\
    \includegraphics[width=0.48\textwidth]{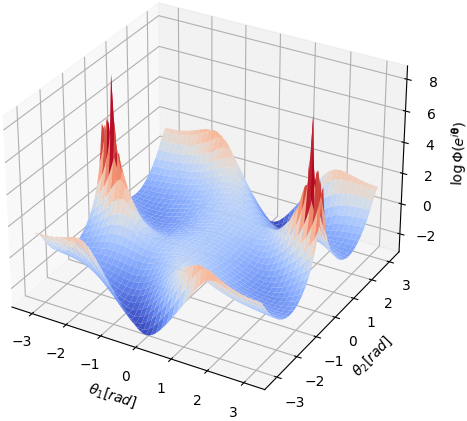}
    \caption{Three sample estimates of the multidimensional spectral density function obtained from the generated WSS process samples using the proposed convex formulation in Algorithm \ref{alg:algorithm1}.}
\label{Fig3}
\end{figure}

Utilizing these finite data records, we compute the multidimensional spectral density estimates via Algorithm \ref{alg:algorithm1}. Three independent sample estimates obtained from Monte Carlo simulations are presented in Figure \ref{Fig3}. The numerical optimization is implemented using the BFGS algorithm provided by the \texttt{scipy.optimize} module, which efficiently converges due to the mathematically proven strict convexity of our formulation. These visualizations provide qualitative evidence of the estimator's ability to consistently recover the spectral shape and analytic smoothness.

For comparative analysis, we establish a baseline utilizing a direct $L_2$-norm minimization strategy, a standard heuristic in applied system identification. Let 
$$
I\left(e^{i \boldsymbol{\theta}}\right)=\sum_{\boldsymbol{k} \in \Omega} \hat{c}_{\boldsymbol{k}} e^{-i(\boldsymbol{k}, \boldsymbol{\theta})}
$$
denote the periodogram of the data. The baseline solves the following constrained optimization problem
$$
\min_{\Lambda \in \mathfrak{L}^{d}_{+}} \left\| \frac{P(e^{i \boldsymbol{\theta}})}{Q(e^{i \boldsymbol{\theta}}, \Lambda)} - I(e^{i \boldsymbol{\theta}}) \right\|_{2}^{2}
$$
where $\Lambda$ is constrained to the set $\mathfrak{L}^{d}_{+}$ to preserve strict positivity. Because the objective function is highly non-convex with respect to $\Lambda$, direct application of the BFGS algorithm frequently converges to sub-optimal local minima.

\begin{figure}[htbp]
    \centering
    \includegraphics[width=0.48\textwidth]{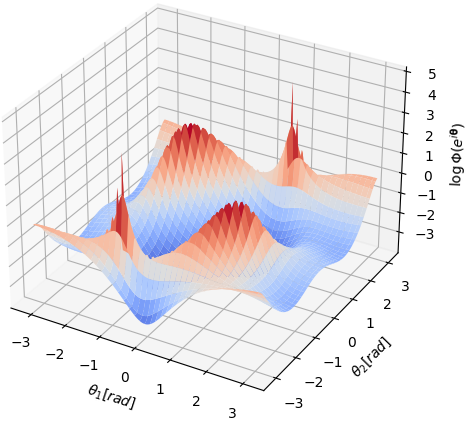}
    \includegraphics[width=0.48\textwidth]{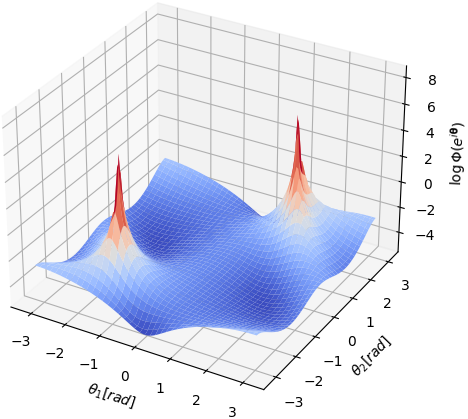}\\
    \includegraphics[width=0.48\textwidth]{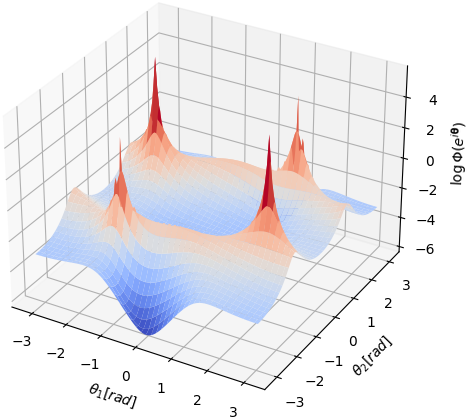}
    \caption{Three sample estimates of the multidimensional spectral density function obtained via the non-convex $L_2$-norm minimization baseline.}
\label{Fig4}
\end{figure}

As shown in Figure \ref{Fig4}, this non-convex baseline yields estimates characterized by severe distortions and non-analytic artifacts. In stark contrast, the algorithm proposed in this paper, which parameterizes the spectrum via trigonometric moments, guarantees existence and uniqueness through a strictly convex framework. 

The simulation results in Figure \ref{Fig3}, corroborated by the theoretical convergence properties established in earlier sections, substantiate the robustness of the proposed estimator. While the current methodology resolves the optimal identification of the denominator $\Lambda$ for a prescribed prior $P(e^{i \boldsymbol{\theta}})$, extending this geometric framework to facilitate the joint convex estimation of both the multidimensional MA parameters and the AR parameters remains a fundamental topic for future research.

\section{A concluding remark}
TMTMP and the corresponding algorithms are an important means of treating the multidimensional spectral density estimation, and multidimensional ARMA modeling for spatiotemporal stochastic processes. However, compared to the trigonometric moment problems with a single dimension, there has not been an explicit condition for the TMTMP to ensure the positiveness of the parametrization of the spectral density in the literature. Moreover, the existence and uniqueness of the solution to the TMTMP are always missing. In this paper, taking real-world applications into consideration, we propose a specific choice of basis functions for the TMTMP and a convex optimization scheme corresponding to it, instead of considering the general TMTMP and analyzing its properties. With this selection of basis functions for the trigonometric moments, the positiveness of the corresponding parametrization of the spectral density function is ensured. The proposed spectral density parametrization is significant, since for any feasible multidimensional trigonometric moment sequence of which the set of basis functions is $\Omega$, there always exists such a parametrization satisfying the trigonometric moment conditions. Based on the proposed parametrization, we propose an estimation scheme for the multidimensional spectral density. The consistency, (asymptotic) unbiasedness, convergence rate and efficiency under a mild assumption are all proved for the spectral estimator, which validates the proposed parametrization and the corresponding estimation scheme. Finally in the numerical simulations, we provide a complete simulation procedure, which include generalizing data samples from the spectral density function and estimating it by the proposed parametrization and the corresponding Algorithm 1. The results are compared to those of estimating the multidimensional spectral density by minimizing the L2-norm of error. The spectral density estimates reveal the advantage of the convex optimization scheme. 

A comprehensive solution to the TMTMP problem from a system and signal processing perspective has been proposed in this paper. However, numerous issues remain to be addressed. In this paper, $\Omega$ is chosen, but it is not the only choice that guarantees the positiveness of the parametrization. A potential research direction is to propose a general characterization for all sets satisfying the positiveness condition. Additionally, as mentioned in the previous section, determining the coefficients of the polynomial in the numerator of the spectral estimate is of great importance.

\bibliographystyle{ieeetr}
\bibliography{Reference}

\begin{IEEEbiography}[{\includegraphics[width=1in,height=1.25in,clip,keepaspectratio]{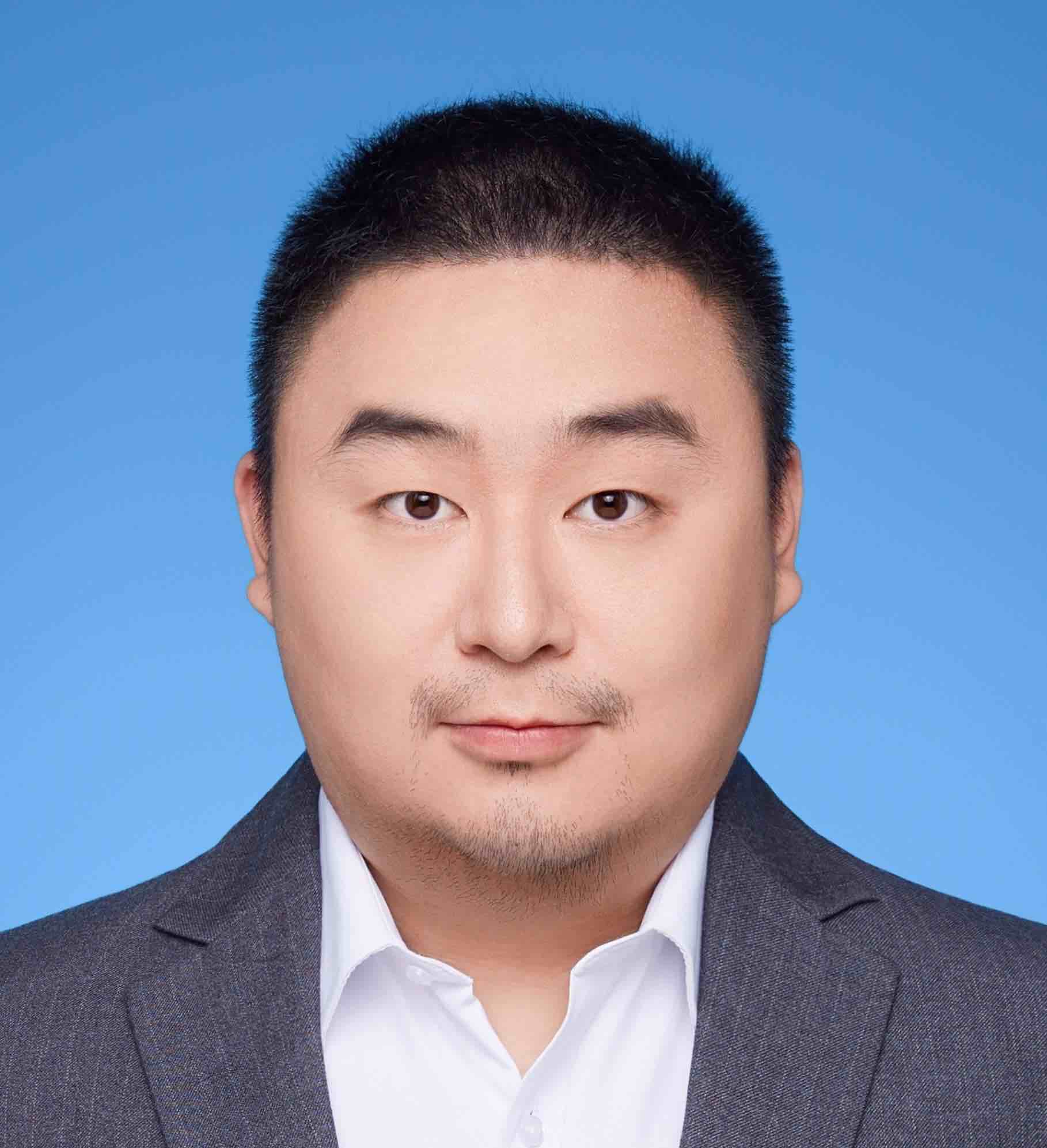}}]{Guangyu Wu} (Member, IEEE) received the Ph.D. degree in Control Science and Engineering (advisor: Anders Lindquist) from Shanghai Jiao Tong University, Shanghai, China, in 2024.

He is an incoming postdoctoral researcher with the Department of Mathematical Sciences, Chalmers University of Technology, Sweden. Prior to that, we worked as a postdoctoral fellow with College of Computing and Data Science, Nanyang Technological University, Singapore. His research interests encompass stochastic control, stochastic filtering, system identification, applied functional analysis, and foundations of machine learning. He serves as an active reviewer for IEEE Transactions on Automatic Control, Automatica, IEEE Control Systems Letters, and the IEEE Conference on Decision and Control. He is a recipient of Eric and Wendy Schmidt AI in Science Postdoctoral Fellowship.
\end{IEEEbiography}

\begin{IEEEbiography}[{\includegraphics[width=1in,height=1.25in,clip,keepaspectratio]{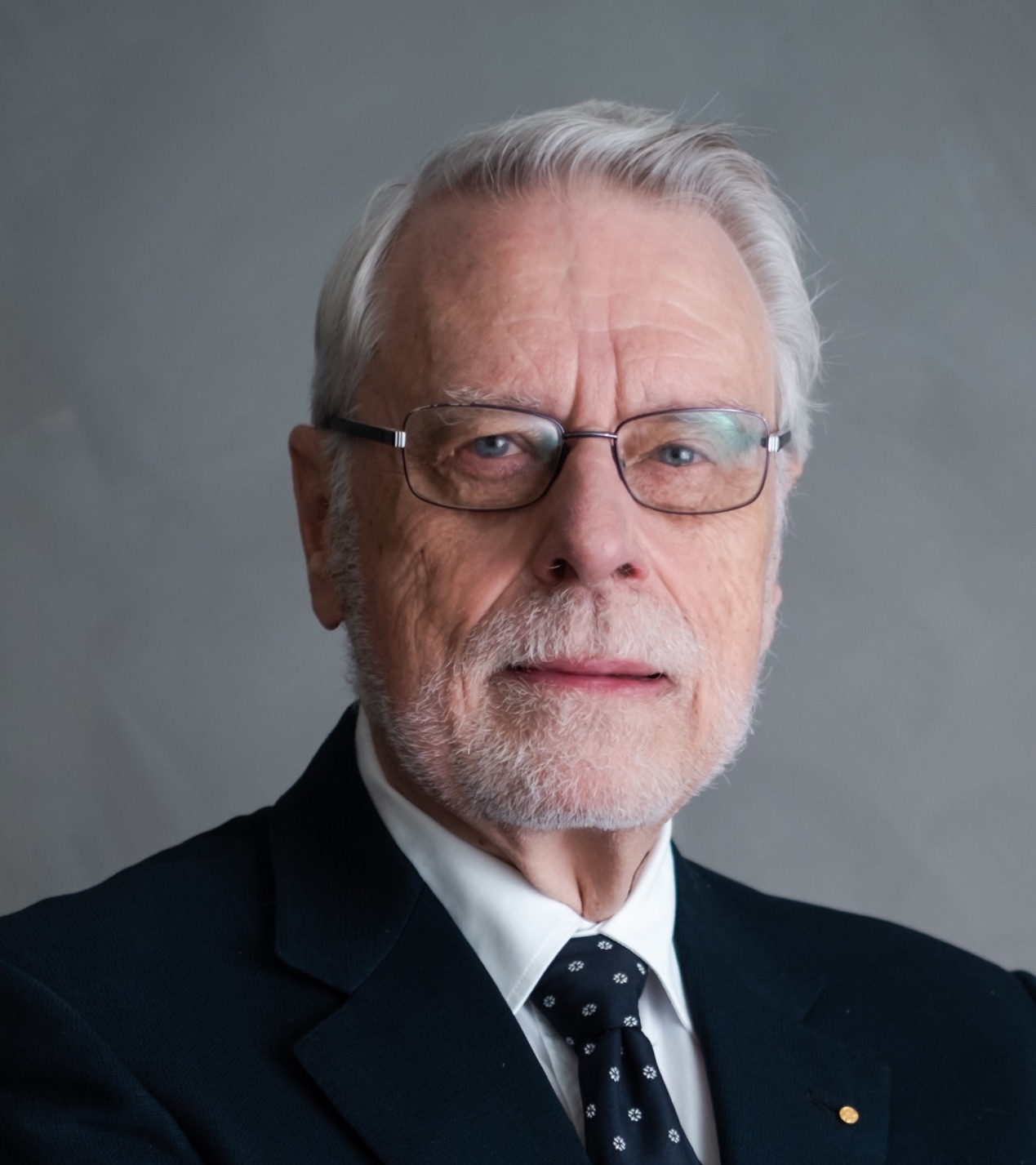}}]{Anders Lindquist} (Life Fellow, IEEE) received the Ph.D. degree in Optimization and Systems Theory from the Royal Institute of Technology (KTH), Stockholm, Sweden, in 1972, an honorary doctorate (Doctor Scientiarum Honoris Causa) from Technion (Israel Institute of Technology) in 2010 and Doctor Jubilaris from KTH in 2022.

He is currently a Distinguished Professor at Anhui University, Hefei, China, Professor Emeritus at Shanghai Jiao Tong University, China, and Professor Emeritus at the Royal Institute of Technology (KTH), Stockholm, Sweden. Before that he had a full academic career in the United States, after which he was appointed to the Chair of Optimization and Systems at KTH.

Dr. Lindquist is a Member of the Royal Swedish Academy of Engineering Sciences, a Foreign Member of the Chinese Academy of Sciences, a Foreign Member of the Russian Academy of Natural Sciences (elected 1997), a Member of Academia Europaea (Academy of Europe), an Honorary Member the Hungarian Operations Research Society, a Life Fellow of IEEE, a Fellow of SIAM, and a Fellow of IFAC. He received the 2003 George S. Axelby Outstanding Paper Award, the 2009 Reid Prize in Mathematics from SIAM, and the 2020 IEEE Control Systems Award, the IEEE field award in Systems and Control.
\end{IEEEbiography}

\end{document}